\documentclass[preprint,12pt]{elsarticle}

\usepackage{amsmath,amssymb}
\usepackage{mathrsfs}
\usepackage{amsthm}
\usepackage{mathtools}
\usepackage{gensymb}
\usepackage{setspace}
\usepackage{extarrows}
\usepackage{geometry}
\usepackage{booktabs}
\usepackage{multirow,multicol}
\usepackage{makecell}
\usepackage{listings}
\usepackage{appendix}
\usepackage{caption,subcaption}
\usepackage{bm}
\usepackage{stack}
\usepackage{subfloat,wrapfig}
\usepackage{xcolor}
\usepackage{paralist}
\usepackage{pdfpages}
\usepackage{indentfirst}
\usepackage[normalem]{ulem}
\usepackage{algorithm}
\usepackage{algorithmic}
\usepackage[hidelinks]{hyperref}

\newtheorem{theorem}{Theorem}[section]

\begin{document}
	
\begin{frontmatter}
	
	\title{A Randomized GMsFEM with Data-Driven Predictors for Parametric Flow Problems in Multiscale Heterogeneous Media}
	
	\author{Wing Tat Leung$^{1}$}
        \ead{wtleung27@cityu.edu.hk}
	\author{Qiuqi Li$^{2}$\corref{cor1}}
	\ead{qli28@hnu.edu.cn}
        \author{Songwei Liu$^{1}$}
        \ead{songwliu2-c@my.cityu.edu.hk}

	    \address{$^1$ Department of Mathematics, City University of Hong Kong, Hong Kong Special Administrative Region}
	\address{$^2$ Department of Mathematics, Hunan University, Changsha 410082, China}

	\cortext[cor1]{Corresponding author}

	\begin{abstract}
		In this paper, we propose a randomized generalized multiscale finite element method (Randomized GMsFEM) for flow problems with parameterized inputs and high-contrast heterogeneous media. The method employs a data-driven predictor to construct multiscale basis functions in two stages: offline and online. In the offline stage, a snapshot space is generated via spectral decompositions, and a reduced matrix is obtained using SVD to predict eigenfunctions. In the online stage, these eigenfunctions are evaluated for new parameter realizations to construct the multiscale space. Furthermore, our approach addresses the complexity of multiple permeability fields with random inputs and multiple multiscale information,  providing accurate and efficient approximations. Moreover, we conduct a rigorous convergence analysis for our Randomized GMsFEM. Finally, we present extensive numerical examples, demonstrating its superior performance compared to the traditional GMsFEM.
	\end{abstract}
	
\end{frontmatter}
	
\section{Introduction}

Multiscale problems are ubiquitous in various physical and engineering applications, such as porous media flow, subsurface transport, and heterogeneous materials. These problems often involve complex systems with multiple source terms and input parameters drawn from high-dimensional spaces. A representative example is flow through heterogeneous porous media, typically modeled by the Darcy equation with uncertain and spatially varying permeability fields. Such permeability fields are often parameterized in a nontrivial manner and exhibit significant uncertainty.
Accurate prediction in these settings requires solving a large number of forward problems under varying source terms and permeability configurations, which can be computationally prohibitive. One promising strategy to alleviate this burden is reduced-order modeling (ROM) \cite{li2017novel,li2020variable,ref17,ref18}, which aims to approximate the input-output map in a lower-dimensional space. In flow problems, the input space typically includes source terms and permeability coefficients, while the output space consists of coarse-grid solutions or integral flow quantities. Since the output often resides on a low-dimensional manifold embedded in a high-dimensional input space, ROM is well-suited for such problems.

Numerous model reduction techniques have been developed, including global, local, and hybrid approaches. The core idea is to identify a low-dimensional approximation space that captures the essential features of the solution manifold. Global methods construct this space from a large number of full-order simulations \cite{ref1,ref2,ref3,ref4}, but they often require recomputation when source terms or boundary conditions change. In contrast, local model reduction methods construct multiscale basis functions in localized (coarse-grid) regions using precomputed offline spaces, avoiding expensive global simulations and enabling more flexible and efficient implementations. 
\cite{hou1997multiscale} constructed the multiscale finite element base functions that are adaptive to the local property of the differential operator. Then it was generalized by Generalized Multiscale Finite Element Method (GMsFEM) \cite{ref9}, which presented a general procedure to construct the offline space and other methods \cite{ref5,ref6,ref8}. GMsFEM brought many advantages. It was thus widely applied in many multiscale problems \cite{chung2014generalized,chung2018multiscale,spiridonov2019generalized} and improved by \cite{chung2018constraint,chung2015residual,chung2018adaptive,li2018multiscale}.
However, significant challenges remain in problems involving multiple permeability fields with random coefficients, especially in the presence of high contrast and strong parameter coupling. These settings introduce complex interactions across scales and parameters, which are difficult to capture using existing local methods.
To address these issues, this paper proposes a randomized GMsFEM that integrates data-driven predictors for constructing parametric multiscale basis functions. The proposed framework leverages statistical learning techniques to approximate local spectral problems under uncertainty, providing a scalable and flexible approach to multiscale problems with strong heterogeneity and parameter coupling.

In multiscale parametric flow problems, solving parametric (or stochastic) eigenvalue problems is essential to capture fine-scale features. In \cite{mach2025solving} two approaches that use the Taylor expansion and Chebyshev expansion were introduced. In \cite{ref23,ref24} a data-driven reduced order model using Gaussian process regression was considered. \cite{lee2024model} proposed framework represents eigenvectors as a low-rank approximation and was implemented using an offline-online strategy. A stochastic collocation method on an anisoptropic sparse grid is proposed in \cite{grubivsic2023stochastic}. All of the above methods have good approximation properties and low computational cost.
However, significant difficulties still remain due to the need to solve the eigenvalue problem of strong coupling between high-contrast information and random parameters. 
To overcome these limitations, we introduce a data-driven model for solving parametric eigenvalue problems within the GMsFEM framework. Specifically, we employ Generalized Polynomial Chaos (gPC) \cite{ref21,ref22} and Gaussian Process Regression (GPR) \cite{williams2006gaussian} to construct predictors that approximate the eigenfunctions associated with new input parameters.

The proposed method operates in two stages: offline and online. In the offline stage, we solve high-fidelity local eigenvalue problems for selected parameters and extract a reduced-order basis using Proper Orthogonal Decomposition (POD). A transformation matrix is then applied to project eigenfunctions into a low-dimensional space, and predictive models are built for the projection coefficients using gPC and GPR. In the online stage, given a new parameter, the trained predictor generates the projection coefficients, which are used to reconstruct the corresponding eigenfunctions by multiplying them with the transformation matrix. These predicted eigenfunctions are then combined with standard multiscale basis functions to form the online basis space.
We refer to the resulting methods as gPC-GMsFEM and GPR-GMsFEM, depending on the choice of predictor. The performance of both approaches is assessed through a series of numerical experiments, assuming input parameters follow standard Gaussian distributions. We consider various settings involving periodic structures and high-contrast media. The results demonstrate that the proposed method achieves high accuracy in approximating fine-scale solutions, with improved efficiency compared to traditional GMsFEM approaches. Moreover, convergence analysis shows that the error decays proportionally to the inverse of the smallest eigenvalue excluded from the reduced space, denoted by $\Lambda_*$.

The remainder of the paper is organized as follows:
In Section~\ref{sec:Preliminaries}, we introduce the necessary preliminaries and notations.
Section~\ref{Multiscale basis based on gPC} presents the construction of the randomized multiscale basis using gPC and GPR.
Numerical examples are provided in Section~\ref{sec:Numerical} to demonstrate the effectiveness of the proposed methods. Conclusions are drawn in Section~\ref{sec:conclusions}.

\section{Preliminaries and notations}
\label{sec:Preliminaries}
In this section, we present some preliminaries and notations for the rest of the paper. Let $L^2(\Omega)$ denote the space of square integrable functions over $\Omega$, and define $H_0^1(\Omega)$. 
Clearly, $L^2(\Omega)$ is equipped with the inner product $(\cdot, \cdot)$, which induces the $L^2$-norm $\|\cdot\|$.

This paper primarily investigates parametric elliptic PDEs, 
\begin{equation}\label{eq:2.1}
    -\nabla\cdot(\kappa(x;\mu)\nabla{u})=f  \ \text{in} \ \Omega,
\end{equation}
where $\kappa(\cdot;\mu)\in L^{\infty}(\Omega)$ is affine-dependent, $\kappa\left(x;\mu\right)=\sum_{q=1}^Q{\Theta}_q\left(\mu\right)\kappa_q\left(x\right) \geq\kappa_{\text{min}}>0$ and $\Theta_q(\mu)$ is the function that only depends on $\mu$. Additionally, the coefficient $\kappa$ incorporates multiscale information, high contrast, and uncertainty. We consider the parameterized elliptic PDEs, which can be formulated as the following weak formulation: find $u\in H^1_0(\Omega)$ such that
\begin{equation}\label{eq:2.2}
    a\left(u,v;\mu\right)=(f,v),\ \ \forall v\in H^1_0(\Omega),
\end{equation}
where $a\left(u ,v ;\mu\right):=\int_{\Omega} \kappa(x;\mu)\nabla u\cdot \nabla v$ is a symmetric bilinear form. $\mu$ is defined in $(\Omega_p, \mathcal{F},P)$, where $\Omega_p \subset \mathbb{R}^{\text{M}}$ and $P$ is arbitrary probability measure. Note that we equip the space $H^1_0(\Omega)$ with the energy norm $\|u\|^2_E = a(u,u; \mu)$. Furthermore, the parametric bilinear form $a\left(\cdot\ ,\cdot\ ;\mu\right)$ is affine with respect to $\mu$. The affine assumption is crucial to enable an offline–online computation decomposition for a multi-query model. When $a\left(\cdot\ ,\cdot\ ;\mu\right)$ is not affine with respect to $\mu$, we can use EIM such that $a\left(\cdot\ ,\cdot\ ;\mu\right)$ can be approximated by an affine representation.  Let $\Xi_{train}$ be a training set, which is a collection of a finite number of samples in $\Omega_p$. Typically the training set is chosen by some randomized sampling methods. We require that the samples in $\Xi_{train}$ are sufficiently scattered in $\Omega_p$. 

\begin{figure}
    \centering
    \includegraphics[width=0.5\linewidth]{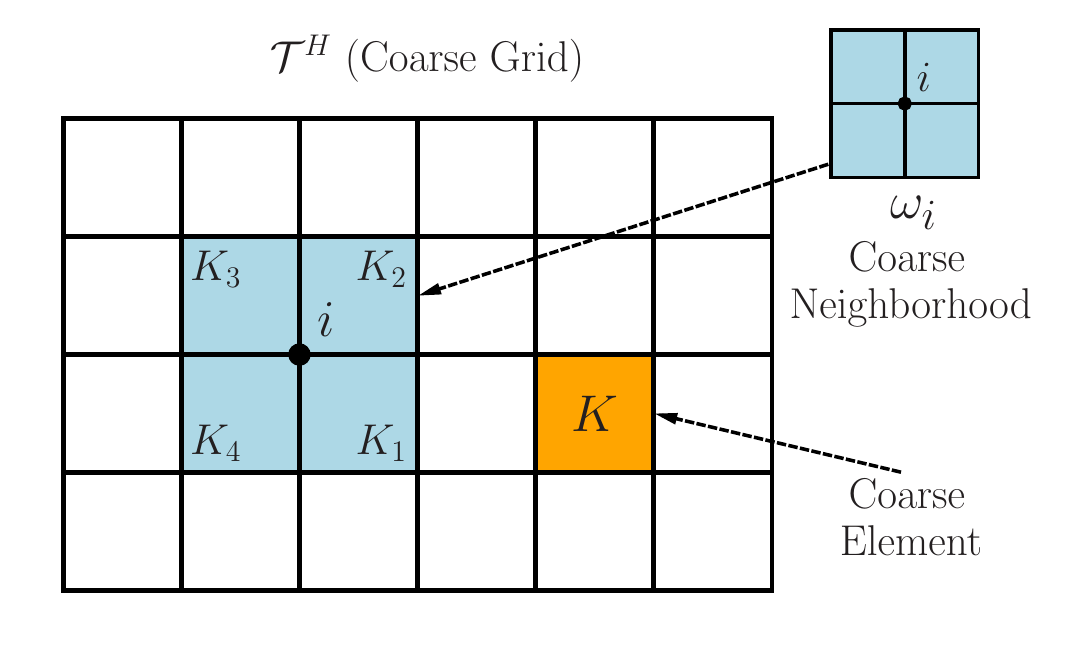}
    \caption{Illustration of a coarse neighborhood and a coarse element. }
    \label{fig:coarse}
\end{figure}

We use $\mathcal{T}^H$ to denote a usual conforming partition of the computational domain $\Omega$. The set $\mathcal{T}^H$ is called the coarse grid and the elements of $\mathcal{T}^H$ are called coarse elements. Moreover, $H>0$ is the coarse mesh size. In this paper, we consider rectangular coarse elements for the ease of discussions and illustrations. The methodology presented can be easily extended to coarse elements with more general geometries. Let $N$ be the number of nodes in the coarse grid $\mathcal{T}^H$, and let $\{x_i|1\leq i \leq N\}$ be the set of nodes in the coarse grid (or coarse nodes for short). For each coarse node $x_i$, we define a coarse neighborhood $\omega_i$ by
\begin{equation}
    \omega_i = \bigcup\{K_j \in \mathcal{T}^H; \ x_i \in \overline{K}_j\}.
\end{equation}
Notice that $\omega_i$ is the union of all coarse elements $K_j \in \mathcal{T}^H$ having the coarse node $x_i$. In Figure~\ref{fig:coarse} we show the definition we describe above.

We let $\mathcal{T}^h$ be a partition of the computational domain $D$ obtained by refining the coarse grid $\mathcal{T}^H$. $\mathcal{T}^h$ is called the fine grid and the fine mesh size is $h$. We remark that the restrictions of the fine grid in coarse neighborhoods will be used to discretize some local problems used to generate local basis functions. Moreover, the fine grid is also used for the computation of a fine scale solution, which is used as a reference solution for comparison purposes. To fix the notations, we will use the standard conforming piecewise linear finite element method for the computation of the fine scale solution. Specifically, we let $V_h$ be the conforming piecewise linear finite space with respect to the fine grid $\mathcal{T}^h$. We will then obtain the fine scale solution $u_h \in V_h$ by solving the following variational problem 
\begin{equation}\label{eq:fine}
    a\left(u_h,v;\mu\right)=(f,v),\ \ \forall v\in V_h.
\end{equation}
Next we present the general idea of GMsFEM. First, we briefly overview the continuous Galerkin formulation of GMsFEM, which has a similar form as the fine scale equation (\ref{eq:fine}). The basis functions are nodal based and have supports on coarse neighborhoods. We denote the $k$-th basis function supported on the coarse neighborhood $\omega_i$ by $\psi_k^{\omega_i}$. For each node $x_i$, we will construct a set of basis functions $\{\psi_k^{\omega_i}|k=1,2,\ldots,l_i\}$, where $l_i$ is the number of basis functions with support in $\omega_i$. These multiscale basis functions are constructed from a local snapshot space and a local spectral decomposition defined on that snapshot space. The snapshot space contains a collection of many basis functions that can be used to capture most of the
fine features of the solution, and the multiscale basis functions $\psi_k^{\omega_i}$ are constructed by selecting the dominant modes of a local spectral problem. Using these multiscale basis functions, the  solution is represented as $u_{ms}(x)=\sum_{i,k}c_k^i\psi_k^{\omega_i}$. Once the basis functions are identified, we can define the approximation space $V_{\text{ms}}$ by the linear span of all basis functions. The GMsFEM solution $u_{\text{ms}}\in V_{\text{ms}}$ can be obtained by the following variational form
\begin{equation}\label{ms}
    a(u_{ms},v; \mu) = (f,v),\ \ \forall v\in V_{\text{ms}}.
\end{equation}
From the above, one sees that the key ingredient of GMsFEM is the construction of local basis functions. We will use the so called offline basis functions, which can be computed in the offline stage.
Moreover, we will construct online basis functions that are parameters dependent and are computed by some predictors we obtained offline. Our results show that  the proposed method has a good approximation of the fine scale solution. Moreover, when the medium has several parameters and high contrasts, our method has a better approximation compared with GMsFEM in the $L^2$ norm.

\subsection{Construction of offline space}
\label{sec:offline}
Next, we will present the construction of the offline basis functions. Let $\omega_i$ be a given coarse neighborhood. We start from a snapshot space $V_{\text{snap}}^{\omega_i}$ which is a set of functions defined on $\omega_i$ and contains all or most necessary components of the fine scale solution restricted to $\omega_i$. A spectral problem is then solved in the snapshot space to obtain the dominant modes in the snapshot space. The snapshot space can be the space of all fine scale basis functions. In general, there are two ways to construct $V_{\text{snap}}^{\omega_i}$. One is the restriction of the conforming space $V_h$ in $\omega_i$ and the resulting basis functions are called spectral basis functions. The other is the set of all $k$-harmonic extensions, and the resulting basis functions are called harmonic basis functions.  

We recall the definition of the snapshot space $V_{\text{snap}}^{\omega_i}$ based on harmonic extensions. Let $J_h(\omega_i)$ be the set of all nodes of the fine mesh $\mathcal{T}^h$ lying on $\partial \omega_i$.
For each fine grid node $x_j \in J_h(\omega_i)$, we define a discrete delta function $\delta_j^h(x)$ defined in $J_h(\omega_i)$ by 
\[
\delta_j^h(x) = \left\{
\begin{array}{cc}
     1, &k=j  \\
     0, &k\neq j 
\end{array}
\right.
, \ \ x_k \in J_h(\omega_i).
\]
The $j$-th snapshot function $\psi_j^{\omega_i,\text{snap}} \in V_{\text{snap}}^{\omega_i}$ for the coarse neighborhood $\omega_i$ is defined as the solution of 
\begin{equation}
\begin{aligned}
    -\nabla\cdot(\kappa(x;\mu)\nabla \psi_j^{\omega_i,\text{snap}}) = 0, \ &\text{in} \ \omega_i, \\
    \psi_j^{\omega_i,\text{snap}} = \delta_j^h(x), \ &\text{on} \ \partial\omega_i.
\end{aligned}
\end{equation}
Obviously, the dimension of $V_{\text{snap}}^{\omega_i}$ is equal to the number of elements in $J_h(\omega_i)$, the set of fine grid nodes lying on $\partial\omega_i$. Let $\{x_j^i\}$ index the elements in $J_h(\omega_i)$, we obtain $V_{\text{snap}}^{\omega_i} = \text{span}\{\psi_j^{\omega_i,\text{snap}}:1\leq j \leq L_i\}$, where $L_i$ is the number of functions we take in each $\omega_i$. 
We then define the corresponding change in the variable matrix $R^{\omega_i}_{\text{snap}}=\left[\Psi_1^{\omega_i,\text{snap}},\ldots,\Psi_{L_i}^{\omega_i,\text{snap}} \right]$, where $\Psi_j^{\omega_i,\text{snap}}$ are considered the columns of the matrix. We observe that employing randomized snapshots alongside oversampling can help decrease the computational expense associated with calculating snapshots. For a detailed explanation, please refer to reference \cite{calo2016randomized}.

We next determine a set of dominant modes from $V_{\text{snap}}^{\omega_i}$, and the resulting lower dimensional space is called the offline space $V_{\text{off}}^{\omega_i}$.
In order to obtain the offline basis functions, we need to perform a reduced space by a series of spectral decompositions. We perform a dimension reduction of the space of snapshots using an auxiliary
spectral decomposition, that is find $(\psi,\lambda) \in V_{\text{snap}}^{\omega_i} \times \mathbb{R}$ such that
\begin{equation}\label{eig}
    \int_{\omega_i}\kappa(x;\mu)\nabla\psi\cdot\nabla\phi = \lambda\int_{\omega_i}\tilde{\kappa}(x;\mu)\psi\phi, \ \forall \phi \in V_{\text{snap}}^{\omega_i},
\end{equation}
where the weight functions $\tilde{\kappa}(x;\mu) = \kappa(x;\mu)H^{-2}$. We arrange the eigenvalues $\lambda_k^{\omega_i}, k =1,2,\ldots,$ obtained by (\ref{eig}) in ascending order. We then select the first $l_i$ eigenfunctions and denote them by $\Psi_1^{\omega_i,\text{off}},\ldots,\Psi_{l_i}^{\omega_i,\text{off}}$. Using these eigenfunctions, we can define 
\[
\phi_k^{\omega_i,\text{off}} = \sum_{j=1}^{L_i}(\Psi_k^{\omega_i,\text{off}})_j\psi_j^{\omega_i,\text{snap}}, \ k=1,2,\ldots,l_i,
\]
where $(\Psi_k^{\omega_i,\text{off}})_j$ denotes the $j$-th component of $\Psi_k^{\omega_i,\text{off}}$. Finally, the offline basis functions for the coarse neighborhood $\omega_i$ are defined by $\psi^{\omega_i,\text{off}}_k = \chi_i\phi^{\omega_i,\text{off}}_k$, where $\chi_i$ is the standard multiscale basis function for the coarse node $x_i$ supported in the neighborhood $\omega_i$ with linear boundary conditions for cell problems. More precisely,
\begin{equation}
\begin{aligned}
    -\nabla\cdot(\kappa(x;\mu)\nabla \chi_i) = 0, \ K \in \omega_i, \\
    \chi_i = g_i, \ \text{on} \ \partial K.
\end{aligned}    
\end{equation}
for all $K\in \omega_i$, where $g_i$ is a continuous function on $\partial K$ and is linear on each edge of $\partial K$. We define the local offline space $V^{\omega_i}_\text{off}$ as the linear span of all $\psi^{\omega_i,\text{off}}_k,k=1,2,\ldots,l_i$.

We remark that one can take $V_{\text{ms}}$ in (\ref{ms}) $V_{\text{off}} = \text{span}\{\psi^{\omega_i,\text{off}}_k|1\leq i \leq N, 1\leq k \leq l_i\}$. We will use offline space as the train space to generate the predictor which will be used to construct the online space.

\section{Randomized GMsFEM}
\label{Multiscale basis based on gPC}

The utilization of $V_{\text{ms}} = V_{\text{off}}$ in (\ref{ms}) presents a promising option in a wide range of scenarios. However, in cases where the problem is parameter dependent, involving multiple multiscale information with random coefficients, one needs to recompute the $V_{\text{ms}}$ when the parameters change. Due to the complexity of the system, recalculating $V_{\text{ms}}$ incurs significant computational costs. In this section, we will propose an offline-online method. We will derive a framework to construct the online basis functions using a predictor when the parameters change. This avoids recomputing the spectral decomposition.

Next, we briefly describe the framework for the construction of multiscale basis functions.
During the offline stage, we employ the POD method to generate snapshots and build the basis of the reduced space. The predictor is utilized to approximate the eigenvalues and projection coefficients in the reduced space, which leads to approximate eigenvectors. The predictor is constructed during the offline stage. In the online stage, given new parameters, we use the predictor to obtain the eigenvalues and projection coefficients of eigenvectors corresponding to the new parameters in the reduced space. The projections are then multiplied by the transformation matrix to obtain the predicted eigenvectors. Finally, the online basis functions are defined by $\psi^{\omega_i,\text{on}} = \chi_i\phi^{\omega_i,\text{on}}$, where $\phi^{\omega_i,\text{on}}$ denotes the predicted eigenvectors.
\subsection{Offline stage}

In the previous sections, we used an auxiliary spectral decomposition (\ref{eig}) to obtain the dominant modes from $V_{\text{snap}}^{\omega_i}$, which is equivalent to the following generalized eigenvalue problem in variational form:
\begin{equation}\label{bilinear}
    a_{\omega_i}(\psi^{\text{snap}},\phi^{\text{snap}}) = \lambda(\mu) s_{\omega_i}(\psi^{\text{snap}},\phi^{\text{snap}}), \ \ \forall \phi^{\text{snap}} \in V^{\omega_i}_{\text{snap}},
\end{equation}
where
\[
a_{\omega_i}(\psi,\phi)=\int_{\omega_i}\kappa(x;\mu)\nabla\psi\cdot\nabla\phi,\ \ 
s_{\omega_i}(\psi,\phi) = \int_{\omega_i}\tilde{\kappa}(x;\mu)\psi\phi.
\]
$a$ and $s$ are symmetric bilinear forms and (\ref{bilinear}) can be written in a matrix form as
\begin{equation}\label{matrix}
    A^{\text{off}}\Psi^{\text{off}}_k = \lambda^{\text{off}}_k S^{\text{off}}\Psi^{\text{off}}_k,
\end{equation}
where
\[
A^{\text{off}}=[a_{mn}]=\int_{\omega_i}\kappa(x;\mu)\nabla\psi^{\text{snap}}_m\cdot\nabla\psi^{\text{snap}}_n=(R^i_{\text{snap}})^TAR^i_{\text{snap}},
\]
and
\[
S^{\text{off}}=[s_{mn}]=\int_{\omega_i}\tilde{\kappa}(x;\mu)\psi^{\text{snap}}_m\psi^{\text{snap}}_n=(R^i_{\text{snap}})^TSR^i_{\text{snap}},
\]
where $A$ and $S$ denote analogous fine scale stiffness and mass matrices defined by 
\[
A_{ij} = \int_{\Omega}\kappa(x;\mu)\nabla\phi_i\cdot\nabla\phi_j, \ \ S_{ij} = \int_\Omega\tilde{\kappa}(x;\mu)\phi_i\phi_j,
\]
where $\phi_i$ is the fine scale basis function for $V_h$. Consider the case where bilinear forms are affine dependent, that is, they can be written as
\begin{subequations}
\begin{gather}
a_{\omega_i}\left(\psi,\phi;\mu\right)=\sum_{q=1}^{\mathrm{\Theta}_a}{\theta_a^q\left(\mu\right)a_{q,\omega_i}\left(\psi,\phi\right)},\ \forall \psi,\phi\in V^{\omega_i}_{\text{snap}},\ \mu\in \Omega_p,\\
s_{\omega_i}\left(\psi,\phi;\mu\right)=\sum_{q=1}^{\mathrm{\Theta}_b}{\theta_b^q\left(\mu\right)b_{q,\omega_i}\left(\psi,\phi\right)},\ \forall \psi,\phi\in V^{\omega_i}_{\text{snap}},\ \mu\in \Omega_p,
\end{gather}
\end{subequations}
where $\theta_a^q(\mu)$ and $\theta_a^q(\mu)$ are the functions that only depend on $\mu$ and here we assume that $\theta_a^q(\mu) = \theta_b^q(\mu) = \Theta_q(\mu) $ as well as $\Theta_a = \Theta_b = Q$.
The bilinear forms on the right-hand side of the equation are parameter independent. That is to say, we can denote the bilinear form as the sum of the products of several parameter dependent functions and parameter independent bilinear forms. The matrix of the corresponding generalized eigenvalue problem can be written as
\begin{subequations}
\begin{gather}
\label{affine}A_{ij} = \int_{\Omega}\kappa(x;\mu)\nabla\phi_i\cdot\nabla\phi_j = \sum_{q=1}^Q{\Theta}_q\left(\mu\right)\int_{\Omega}\kappa_q\left(x\right)\nabla\phi_i\cdot\nabla\phi_j,\\
S_{ij} = \int_\Omega\tilde{\kappa}(x;\mu)\phi_i\phi_j = \sum_{q=1}^Q{\Theta}_q\left(\mu\right)\int_{\Omega}\tilde{\kappa}_q(x)\phi_i\phi_j.
\end{gather}
\end{subequations}

\subsubsection{Construction of reduced space}\label{Dimension reduction}
Previously, we defined $\mathrm{\Xi}_{train}$ as a collection of a finite number of samples in $\Omega_p$, which can be written as $\{\mu_j :j = 1,2,\ldots,n_s\}$. Sample points can be selected using various strategies. For example, we can use uniform tensor sampling, Latin hypercube sampling, or Smolyak sparse sampling. 
For each $\mu_j$, we construct the corresponding snapshot space, solving the generalized eigenvalue problem (\ref{matrix}) and obtain the $l_i$ eigenpairs $(\lambda_i^{\omega_i},\phi_i^{\omega_i})\in \mathbb{R} \times V_{\text{snap}}^{\omega_i}$ for each $\omega_i$. Let $S_n$ be the snapshot matrix of size $n_v$, which contains the eigenvectors at the sample points $\mu_1,\mu_2,\ldots,\mu_{n_s}$, that is
\[
S_n=
\left[\ldots|s^{i,j}_k|\ldots\right]\in \mathbb{R}^{n_h\times n_v}, \ s_k^{i,j}
=
\left(\Psi_{k}^{\omega_i,\mu_j}\right)^T, 
\]
where $n_h = L_i$ and $n_v = l_i \times N \times n_s$, $N$ is the number of coarse elements defined above. Obviously, we can define different $l_i$ in each $\omega_i$. In this case, we assume that we take the same number of eigenvectors. Next, we want a low dimensional matrix of $S_n$. There are two commonly used methods for dimensionality reduction, namely the POD method and the greedy method. In this paper, we use POD to obtain a reduced matrix, namely $S_N$. We denote $S_n$ in SVD form as follows:
\begin{subequations}
\begin{gather}
S_n=U\mathrm{\Sigma}Z^T, ~U=\left\{\beta_1,\ldots,\beta_{n_h}\right\}\in R^{n_h\times n_h},\\
Z=\left\{\delta_1,\ldots,\delta_{n_v}\right\}\in R^{n_v\times n_v},\ \mathrm{\Sigma}=diag\left\{\sigma_1,\ldots\sigma_r,0,..,0\right\}\in R^{n_h\times n_v}, \\
\sigma_1\geq\sigma_2\geq\ldots\geq\sigma_r>0.
\end{gather}
\end{subequations}
where $r$ is the rank of matrix $S_n$. The $N_h$ left eigenvectors of $U$ will be the basis of $V$, i.e. ${\beta_1,\ldots,\beta_{N_h}}$ will be the basis of the low dimensional subspace $V$. We generally do not directly use SVD to find the left singular vectors unless $n_v\geq n_h$. Next, we will describe the calculation process of the POD basis. We have  
\begin{subequations}
\begin{gather}
S_n\delta_i=\sigma_i\beta_i,\ S_n^T\beta_i=\sigma_i\delta_i,\ i=1,2,\ldots,r,\\
S_n^TS_n\delta_i=\sigma_i^2\delta_i,\ S_nS_n^T\beta_i=\sigma_i^2\beta_i,\ i=1,2,\ldots,r.
\end{gather}
\end{subequations}
Here, the matrix $C=S_n^TS_n$ is called the correlation matrix. The POD basis $V$ with dimension $N_h\le n_s$ is defined as the set of the first $N_h$ left singular vectors ${\beta_1,\ldots,\beta_{N_h}}$ of $S_n$, or the equivalent set of vectors as follows
\begin{equation}
    \beta_j=\frac{1}{\sigma_j}S_n\delta_j,\ 1\le j\le N_h,
\end{equation}
which are obtained from the first $N_h$ eigenvectors of the correlation matrix $C$, denoted as $\{\delta_1,\ldots,\delta_{N_h}\}$. Therefore, we do not use the SVD of matrix $S_n$, and we can solve the eigenvalue problem $C\delta=\sigma^2\delta$. Note that $n_s$ is usually much smaller than $n_h$.

Following the usual practice, we choose the POD dimension $N_h$ as the smallest positive integer that satisfies the following conditions
\begin{equation}
I(N_h)=\frac{\sum_{i=1}^{N_h}\sigma_i^2}{\sum_{i=1}^r\sigma_i^2}\geq1-\epsilon.
\end{equation}
When we obtain the POD basis $V$ with dimension $N_h$, we can calculate the reduced matrix $S_N$.
The relationship between the snapshot matrix $S_n$ and the reduced order snapshot matrix $S_N$ is as follows
\begin{equation}
    S_N=V^TS_n=[V^Ts^{1,1}_1|\ldots|V^Ts_{l_i}^{N,n_s}]\in R^{N_h \times n_v}.
\end{equation}

\subsubsection{Construction of predictor }
\label{data-driven model}
In this section, we will describe the data-driven predictor used to compute the eigenvalue and the eigenvectors for each $\omega_i$. This method can be extended to multiple eigenpairs by considering relevant feature vectors in the snapshot matrix. We will also address the challenge of simultaneously searching for multiple feature vectors.

Firstly, use some sampling techniques to select sample parameters and form the training set, that is $\mathrm{\Xi}_{train} =\{ {\mu}_1,\mu_2,\ldots,\mu_{n_s}\}$. For each $\omega_i$ and $\mu_j$, We calculate the $l_i$ eigenpairs of the high fidelity eigenvalue problem in the snapshot space $V^{\omega_i}_{\text{snap}}$ of $L_i$ dimensions in ascending order for each $\omega_i$ and use them to construct snapshot matrix $S_n$. Then, applying the POD method to the snapshot matrix $S_n$ to form the reduced matrix $S_N$.

Generalized polynomial chaos (gPC) is a mathematical tool used for uncertainty quantification (UQ). It represents a random variable as a linear combination of orthogonal polynomials. The selection of orthogonal polynomials usually depends on the distribution type of random variables.  For Gaussian-distributed random variables, Hermite polynomials are used, while Legendre polynomials are suitable for uniform distributions.  gPC flexibly handles non-Gaussian random variables more flexibly through different polynomial families.

Given $\mu=(\mu^1,\ldots,\mu^M)\in \Omega_p \subset \mathbb{R}^{\text{M}}$, the multivariate orthogonal polynomials are defined as:
\[
      \eta_{\bm\alpha}(\mu) = \prod\limits_{m=1}^M\eta_{\alpha_m}^{(m)}(\mu^m),\ \ \ \ \ \bm\alpha =(\alpha_1,\alpha_2,\ldots,\alpha_m)\in\mathbb{N}_0^{\text{M}},
\]
where $\bm\alpha$ is a multi-index and $\mathbb{N}_0^{\text{M}}$ is  a multi-index set of countable cardinality. For numerical computations, we truncate the polynomial set to a finite order $p$, such that  $|\bm\alpha|= \sum_{j=1}^M\alpha_j \leq p$. For a random eigenfunction $\phi_j^{\omega_i}(\mu)$, the  gPC expansion can be obtained as
\begin{equation}\label{eq:4.1}
\phi_j^{\omega_i}(\mu) = \sum\limits_{\bm\alpha\in\mathbb{N}_0^{\text{M}}}c^{\omega_i,j}_{\bm\alpha}\eta_{\bm\alpha}{(\mu)},
\end{equation}
where $c^{\omega_i,j}_{\bm\alpha}$ is the expansion coefficient to be solved, and $\eta_{\bm\alpha}(\mu)$ is the orthogonal polynomial corresponding to the random variable $\mu$. 
The gPC coefficients are computed using the discrete least squares method. 
Considering the following reduced matrix, we obtain the eigenvectors
\begin{equation}   
S_N=[V^Ts^{1,1}_1|\ldots|V^Ts^{1,n_s}_1|\ldots|V^Ts_{l_i}^{N,n_s}]\in R^{N_h \times n_v}.
\end{equation}
We extract the corresponding training data of each eigenvector from each row. Let $P_{k,N}^{\omega_i}:\Omega_p \rightarrow R^{N_h}$ be a mapping defined by $\mu\rightarrow \xi^{\omega_i,\text{on}}_{k,N}(\mu)$, where $\xi^{\omega_i,\text{on}}_{k,N}(\mu)$ is the online eigenvector we want to approximate. Using the training data 
\[
\left\{\left(\mu_j,{\Psi_{k,N}^{\omega_i}}\left(\mu_j\right)\right):\ j=1,\ldots,n_s\right\},
\]
we wish to approximate this mapping $P_{k,N}^{\omega_i}$ by a gPC function $\widehat{P}_{k,N}^{\omega_i}:\Omega_p \rightarrow \widehat{\xi}_{k,N}^{\omega_i,\text{on}}$. In the offline stage, a gPC-based prediction function  is constructed. In the online stage, for any given parameter $\mu_*$, all gPC functions $\widehat{P}_{k,N}^{\omega_i}$ are evaluated. We first compute the reduced eigenvector and then multiplying the transform matrix $V$ to obtain the high fidelity eigenvector $\widehat{\xi}_{k}^{\omega_i,\text{on}}$, which is the corresponding coefficient vector in terms of the basis for snapshot space $V^{\omega_i}_{\text{snap}}$.
The predicted coefficients form a predicted eigenfunction
\[
\widehat{\phi}_k^{\omega_i,\text{on}} = \sum_{j=1}^{L_i}(\widehat{\xi}_k^{\omega_i,\text{on}})_j\psi_j^{\omega_i,\text{snap}}, \ k=1,2,\ldots,l_i,
\]
where $(\widehat{\xi}_k^{\omega_i,\text{on}})_j$ denotes the $j$-th component of $\widehat{\xi}_k^{\omega_i,\text{on}}$. We can use this vector as a reference to extract all eigenvectors and construct corresponding eigenfunctions.

At the end of next part, we will give an algorithm to construct the offline stage of randomized GMsFEM. We did not specify a specific predictor because one can use different predictors due to the problem and parameters. In Appendix~\ref{GPR}, we demonstrate how to construct the predictor via GPR. Moreover, a predictor that does not appear in the article can be chosen carefully.

\subsection{Online stage}
In this section, we discuss online multiscale spaces that are constructed for parameter-dependent problems.
Given a new parameter $\mu_*$, we use the gPC predictors $\widehat{P}_{k,N}^{\omega_i}$, multiplied by the transform matrix $V$, to generate the high fidelity predicted eigenvectors $\widehat{\xi}_{k}^{\omega_i,\text{on}}$. The predicted eigenfunctions are given by  $\widehat{\phi}_k^{\omega_i,\text{on}}(\mu) = \sum_{j=1}^{L_i}(\widehat{\xi}_k^{\omega_i,\text{on}})_j\psi_j^{\omega_i,\text{snap}}, \ k=1,2,\ldots,l_i$ for each $\omega_i$. The online basis functions for the coarse neighborhood $\omega_i$ are defined by $\psi^{\omega_i,\text{on}}_k(\mu) = \chi_i\widehat{\phi}^{\omega_i,\text{on}}_k(\mu)$, where $\chi_i$ is the standard multiscale basis function for the coarse node $x_i$ supported in the neighborhood $\omega_i$ with linear boundary conditions for cell problems. 
We then define the online space as $V_{\text{on}} = \text{span}\{\psi^{\omega_i,\text{on}}_k:1\leq i \leq N, 1 \leq k \leq l_i\}$. For any $ u_H \in V_{\text{on}}$, we have $u_H = \sum_{i=1}^N\sum_{k=1}^{l_i}c_i^k\psi_k^{\omega_i,\text{on}}$. To simplify the symbols, we denote $u_H = \sum_{i=1}^{N_v}c_i\psi_i^{\text{on}}$, where $N_v$ is the number of basis functions. 
The coarse solution $u_{H}\in V_{\text{on}}$ can be obtained by the following variational form
\begin{equation}\label{on}
    a(u_{H},v) = (f,v),\ \ \forall v\in V_{\text{on}}.
\end{equation}

For each snapshot basis function $\psi_j^{\omega_i,\text{snap}}$, let $\Psi_j^{\omega_i,\text{snap}}$ be the corresponding coefficient vector in terms of the basis for fine grid space. We can put all snapshot functions using a matrix representation $R_{\text{snap}}^{\omega_i} = \left[\Psi_1^{\omega_i,\text{snap}},\ldots,\Psi_{L_i}^{\omega_i,\text{snap}} \right]$. Similarly, we put all online basis functions in matrix form, that is $C^{\omega_i}_{\text{on}}= \left[\zeta_1^{\omega_i,\text{on}},\ldots,\zeta_{l_i}^{\omega_i,\text{on}} \right]$. The matrix form of (\ref{on}) can be written as
\begin{equation}\label{matrixon}
    C_{\text{on}}^TR^T_{\text{snap}}AR_{\text{snap}} C_{\text{on}}u_H^{\text{discrete}}= C_{\text{on}}^TR^T_{\text{snap}}F,
\end{equation}
where $u_H^{\text{discrete}}$ is the discrete version of $u_H$ and $A$ denotes fine scale stiffness matrix, $F$ is the vector form of source term. According to the affine-dependence of $\kappa(\cdot;\mu)$, we express $A$ as 
\[
A_{ij} = \int_{\Omega}\kappa(x;\mu)\nabla\phi_i\cdot\nabla\phi_j = \sum_{q=1}^Q{\Theta}_q\left(\mu\right)\int_{\Omega}\kappa_q\left(x\right)\nabla\phi_i\cdot\nabla\phi_j, \ \forall \phi_i \in V_h.
\]
Thus (\ref{matrixon}) can be expressed as 
\begin{equation}\label{matrixmu}
    \left(\sum_{q=1}^Q \Theta_q(\mu) C_{\text{on}}^T\widehat{A_q}C_{\text{on}}\right) u_H^{\text{discrete}}= C_{\text{on}}^T\widehat{F},
\end{equation}
where $\widehat{A_q} = R^T_{\text{snap}}A_qR_{\text{snap}}$, $\widehat{F} = R^T_{\text{snap}}F$ and we denote $A_q$ by
$
A_q =\int_{\Omega}\kappa_q\left(x\right)\nabla\phi_i\cdot\nabla\phi_j$.

The benefit of doing so is that we can calculate and store the reduced matrix $\widehat{A_q}$ in advance in the offline stage for $q=1,2,\ldots,Q$. In the online stage, when the parameter changes, we only need to recalculate the online basis functions and obtain $C_{\text{on}}$. We also need to recalculate the parameter terms $\Theta_q(\mu)$ for $q=1,2,\ldots,Q$. Then we couple the matrices together and obtain the coarse grid solution $u_H$.

In addition, we can adopt a different strategy to solve (\ref{on}). Given a new parameter $\mu_*$, we use predictors $P_{k,N}^{\omega_i}$ to obtain the low fidelity predicted eigenvectors $\widehat{\xi}_{k,N}^{\omega_i,\text{on}}$, which are the corresponding coefficients vectors in terms of the basis for POD space defined by $V_{\text{POD}} = \text{span}\{\psi_k^{\text{p}}:1\leq k \leq N_h \}$. 
The simulation of low fidelity predicted eigenfunctions can be written as $
\widehat{\phi}_{k,N}^{\omega_i,\text{on}} = \sum_{j=1}^{N_h}(\widehat{\xi}_{k,N}^{\omega_i,\text{on}})_j\psi_j^{\text{p}}, \ k=1,2,\ldots,l_i \ .$
We then define the online space as $V_{\text{on}} = \text{span}\{\psi^{\omega_i,\text{on}}_{k,N} = \chi_i\widehat{\phi}^{\omega_i,\text{on}}_{k,N} :1\leq i \leq N, 1 \leq k \leq l_i\}$. And we denote $\psi^{\omega_i,\text{on}}_{k,N}(\mu)$ for each $\omega_i$ by
$
\psi^{\omega_i,\text{on}}_{k,N}(\mu) = \sum_{j=1}^{N_h}(\zeta_{k,N}^{\omega_i,\text{on}})_j\psi_j^{\text{p}}, \ k=1,2,\ldots,l_i
$. Similarly, we define $C^{\omega_i}_{\text{on}}= \left[\zeta_{1,N}^{\omega_i,\text{on}},\ldots,\zeta_{l_i,N}^{\omega_i,\text{on}} \right]$ 
and 
$R_{\text{p}} = \left[ \Psi_1^{\text{p}},\ldots,\Psi_{N_h}^{\text{p}}  \right]$.
The matrix form of (\ref{on}) can be expressed as 
\begin{equation}\label{matrixpod}
        C_{\text{on}}^TR_{\text{p}}^TR^T_{\text{snap}}AR_{\text{snap}}R_{\text{p}} C_{\text{on}}u_H^{\text{discrete}}= C_{\text{on}}^TR^T_{\text{p}}R^T_{\text{snap}}F.
\end{equation}
As the presentation above, we can express the matrix form as (\ref{matrixmu}), that is
\begin{equation}\label{mumatrixpod}
    \left(\sum_{q=1}^Q \Theta_q(\mu) C_{\text{on}}^T\widehat{A_q}C_{\text{on}}\right) u_H^{\text{discrete}}= C_{\text{on}}^T\widehat{F},
\end{equation}
where $\widehat{A_q} = {R}_{\text{p}}^T{R}_{\text{snap}}^TA_{q}{R}_{\text{snap}}{R}_{\text{p}}$ and $\widehat{F} = R^T_{\text{p}}R^T_{\text{snap}}F$.
We can calculate and store a small matrix $\widehat{A}$ in advance in the offline stage so that we avoid recomputing the matrix when parameter changes. Moreover, a small matrix can improve the computing speed and reduce computational costs.

Finally, we present the convergence theorem of the method we proposed above. For convenience, we name our method gPC-GMsFEM. The proof and details of the theorem are given in subsection~\ref{convergence}.

\begin{theorem}\label{thm}
Let u and $u_H$ be the solutions to problems (\ref{eq:2.2}) and (\ref{on}), respectively. There holds 
\begin{equation}
E\left(\int_\Omega\kappa|\nabla(u-u_H)|^2\right) \leq \frac{C_1}{\Lambda_*} + C_2H^2 +
C_3 \widetilde{\kappa}_{\text{max}}\gamma + C_3 \widetilde{\kappa}_{\text{max}}|\delta| + C_4\widetilde{\kappa}_{\text{max}}\sum_{i,j,k}\int_P(\epsilon_N^{\omega_i,j,k})^2dP(\mu),
\end{equation}
where $E$ is the expectation, $\Lambda_*$ is the smallest eigenvalue in $\Omega_p$ and not included in the basis space and $\gamma$ satisfies $\sum_{k=1}^r\lambda_k/(\sum_{k=1}^{\infty}\lambda_k) \geq 1- \gamma$, where r is the terms of POD expansion and $\lambda_k$ denotes the eigenvalue of POD expansion. $\delta$ is a small number indicating the error of sampling.   $\epsilon_N^{\omega_i,j,k}$ is the  error between the POD expansion coefficients and the gPC functions expanded by N terms, $P$ is the probability measure.
\end{theorem}

At the end of this part, we give an algorithm to construct both the offline and the online stage of randomized GMsFEM. We will calculate the predictors $\widehat{P}_{k,N}^{\omega_i}$ as well as reduced matrix $\widehat{A_q}$ in the offline stage. Once the new parameter, source term and boundary condition are given, the online space directly generated by the predictor is used to solve the weak formulation. 
\begin{algorithm}
\caption{Offline stage of Randomized GMsFEM}
\label{alg:example}
\begin{algorithmic}[1] 
    \renewcommand{\algorithmicrequire}{\textbf{Input:}}
    \REQUIRE Randomly sampled parameters ${\mu}_1,\mu_2,\ldots,\mu_{n_s}$  and corresponding $\kappa(x;\mu_i)$

    \STATE Construct fine grid $\mathcal{T}^h$ and coarse grid $\mathcal{T}^H$
    \STATE Construct snapshot matrix $R_{\text{snap}}$  
    \STATE Compute the $l_i$ eigenvectors for each $\omega_i$ and $\mu_j$ to form the snapshot matrix $S$
    \STATE Construct the reduced matrix $S_N = V^TS$ and the POD matrix $R_{\text{p}}$
    \STATE Construct predicted function $\widehat{P}_{k,N}^{\omega_i}$ and the small matrix $\widehat{A_q}$

   \RETURN $\widehat{P}_{k,N}^{\omega_i}$, $\widehat{A_q}$, $R_{\text{snap}}$, $R_{\text{p}}$
\end{algorithmic}
\end{algorithm}

\begin{algorithm}
\caption{Online stage of Randomized GMsFEM}
\label{alg}
\begin{algorithmic}[1] 
    \renewcommand{\algorithmicrequire}{\textbf{Input:}}
    \REQUIRE New parameter $\mu_*$, $\widehat{P}_{k,N}^{\omega_i}$, $\widehat{A_q}$, $R_{\text{snap}}$, $R_{\text{p}}$, $F$ and the boundary condition
    \STATE Compute the function $\Theta_q(\mu_*)$
    \STATE Construct the online matrix $C_{\text{on}}$ and online vector $\widehat{F}$
    \STATE Solve the coarse grid problem (\ref{mumatrixpod}) for given source term and boundary condition
    
    \RETURN $u_H^{\text{discrete}}$
\end{algorithmic}
\end{algorithm}

\subsection{The convergence study of gPC-GMsFEM}\label{convergence}
In this section, we study the convergence of gPC-GMsFEM. We will discuss the expectation of the energy error. Assume that the parameters $\mu$ follow a certain multidimensional probability measure $P(\mu)$ and are independent of each other. Let $I^{\omega_i}u$ denote the first $l_i$ terms of spectral expansion of $u$ in terms of eigenfunctions of $-div(\kappa\nabla\phi_j^{\omega_i})= H^{-2}\lambda_j\kappa\phi_j^{\omega_i}$ in $\omega_i$.
Next, we take $\hat{\phi}_j^{\omega_i}$ to be the $r\text{th}$ order POD expansion of $\phi_j^{\omega_i}$ in each $\omega_i$, that is $\hat{\phi}_j^{\omega_i} = \sum_{k=1}^r a_k^{{\omega_i},j}(\mu)\psi_k(x)$, where $\psi_k(x)$ is defined by $\int_{\Omega}C(x,x')\psi_k(x')dx' = \lambda_k\psi_k(x)$. Then we define $\hat{\phi}_{j,1}^{\omega_i}$ as the gPC approximation function and $\hat{\phi}_{j,2}^{\omega_i}$ as the predicted function.
Let $I_0^{\omega_i}u$ denote the $l_i$ terms of combination of $\hat{\phi}_j^{\omega_i}$ in $\omega_i$. 
Similarly, we can define $I_1^{\omega_i}u$ as well as $I_2^{\omega_i}u$. 
Finally we define $\chi_i$ as a partition of unity function satisfying $\sum_i \chi_i = 1$, where $|\nabla\chi_i|\leq \frac{1}{H}$.
We now proceed to prove (\ref{thm}). 

\begin{proof}
Using the Galerkin finite element formulation on a coarse grid, we have
\begin{equation}
\scalebox{0.8}{$\displaystyle
\begin{split}
    \int_\Omega\kappa|\nabla(u-u_H)|^2 &\leq \sum_i\int_{\omega_i}\kappa|\nabla(\chi_i(u-I_2^{\omega_i}u))|^2\\
    &\leq \sum_i \int_{\omega_i}\kappa|\nabla(\chi_i)|^2(u-I_2^{\omega_i}u)^2 + \sum_i \int_{\omega_i}\kappa\chi_i^2|\nabla(u-I_2^{\omega_i}u)|^2.
\end{split}$}
\label{eq:B.1}
\end{equation}
Next we estimate the second term on the right hand side of (\ref{eq:B.1}). We denote in each $\omega_i$
\begin{equation}
    -\text{div}(\kappa\nabla(u-I_2^{\omega_i}u))=g,
\label{eq:B.2}
\end{equation}
where $g$ is the residual in the approximation. Multiplying both sides of (\ref{eq:B.2}) by $\chi_i^2(u-I_2^{\omega_i}u)$ and taking integration, we have
\begin{equation}
\scalebox{0.8}{$\displaystyle
\begin{split}
    \int_{\omega_i}\kappa\chi_i^2|\nabla(u-I_2^{\omega_i}u)|^2 + 2\int_{\omega_i}\kappa\chi_i\nabla\chi_i(\nabla(u-I_2^{\omega_i}u))(u-I_2^{\omega_i}u) = \int_{\omega_i}g\chi_i^2(u-I_2^{\omega_i}u).
\end{split}$}
\end{equation}
From this inequality, we deduce
\begin{equation}
\scalebox{0.8}{$\displaystyle
\begin{split}
\int_{\omega_i}\kappa\chi_i^2|\nabla(u-I_2^{\omega_i}u)|^2 &= -2\int_{\omega_i}\kappa\chi_i\nabla\chi_i(\nabla(u-I_2^{\omega_i}u))(u-I_2^{\omega_i}u) +\int_{\omega_i}g\chi_i^2(u-I_2^{\omega_i}u) 
\\&\leq 2\int_{\omega_i}\kappa|\chi_i||\nabla\chi_i||\nabla(u-I_2^{\omega_i}u)||(u-I_2^{\omega_i}u)| + |\int_{\omega_i}g\chi_i^2(u-I_2^{\omega_i}u) |
\\&\leq  C \left( \int_{\omega_i}\kappa|\nabla\chi_i|^2(u-I_2^{\omega_i}u)^2 + \int_{\omega_i}\kappa|\chi_i|^2|\nabla(u-I_2^{\omega_i}u)|^2 \right) + |\int_{\omega_i}g\chi_i^2(u-I_2^{\omega_i}u)|
\\&\leq \hat{C} \left(\int_{\omega_i}\kappa|\nabla\chi_i|^2(u-I_2^{\omega_i}u)^2+ |\int_{\omega_i}g\chi_i^2(u-I_2^{\omega_i}u)| \right).
\end{split}$}
\end{equation}
Then we have
\begin{equation}
\scalebox{0.8}{$\displaystyle
\begin{split}
    \int_\Omega\kappa|\nabla(u-u_H)|^2 &\leq (1+\hat{C})\sum_i\int_{\omega_i}\kappa|\nabla\chi_i|^2(u-I_2^{\omega_i}u)^2 + \hat{C}\sum_i|\int_{\omega_i}g\chi_i^2(u-I_2^{\omega_i}u)|
    \\ &\leq C \left( \sum_i\int_{\omega_i}\kappa|\nabla\chi_i|^2(u-I_2^{\omega_i}u)^2 + 
    \sum_i\int_{\omega_i}(\kappa|\nabla\chi_i|^2)^{-1}g^2 \right),
\end{split}$}
\end{equation}
where the last step is derived by $ab \leq a^2 + b ^2$ and $\chi_i \leq 1$.
Then we have
\begin{equation}
\scalebox{0.77}{$\displaystyle
\begin{split}
    \int_\Omega\kappa|\nabla(u-u_H)|^2 &\leq C\sum_i\int_{\omega_i}\kappa|\nabla\chi_i|^2(u-I_2^{\omega_i}u)^2 + C\sum_i\int_{\omega_i}(\kappa|\nabla\chi_i|^2)^{-1}g^2
    \\ &\leq C \sum_i\left(\int_{\omega_i}\kappa|\nabla\chi_i|^2(u-I^{\omega_i}u)^2 +\int_{\omega_i}\kappa|\nabla\chi_i|^2(I^{\omega_i}u-I_0^{\omega_i}u)^2 +\int_{\omega_i}\kappa|\nabla\chi_i|^2(I_0^{\omega_i}u-I_1^{\omega_i}u)^2 \right) \\& \ \ \ 
    +C \sum_i\int_{\omega_i}\kappa|\nabla\chi_i|^2(I_1^{\omega_i}u-I_2^{\omega_i}u)^2
    + C\sum_i\int_{\omega_i}(\kappa|\nabla\chi_i|^2)^{-1}g^2
    \\ &\leq C \sum_i\frac{1}{\lambda_{l_i+1}^{\omega_i}}\int_{\omega_i}\kappa|\nabla(u-I^{\omega_i}u)|^2
    + C \sum_i\int_{\omega_i}\kappa|\nabla\chi_i|^2(I^{\omega_i}u-I_0^{\omega_i}u)^2
    \\& \ \ \ 
    + C\sum_i\int_{\omega_i}\kappa|\nabla\chi_i|^2(I_0^{\omega_i}u-I_1^{\omega_i}u)^2
    + C\sum_i\int_{\omega_i}\kappa|\nabla\chi_i|^2(I_1^{\omega_i}u-I_2^{\omega_i}u)^2
    \\& \ \ \ 
    + C\sum_i\int_{\omega_i}(\kappa|\nabla\chi_i|^2)^{-1}g^2 .
\end{split}$}
\end{equation}
We assume that there exist a global function $G$ satisfying $\int_{\Omega}G^2 \leq C$ such that $\int_{\omega_i}(H^{2}\kappa|\nabla\chi_i|^2)^{-1}g^2 \leq \int_{\omega_i}G^2 $. Note that $\int_{\omega_i}\kappa|\nabla(u-I^{\omega_i}u)|^2 \leq \int_{\omega_i}\kappa|\nabla u|^2$. Taking $\Lambda_* = \text{min}_{\omega_i, \mu}\lambda_{l_i+1}^{\omega_i}(\mu)$ and assuming that $\int_{\Omega}G^2$ and $\int_{\Omega}\kappa|\nabla u|^2$ are bounded, we have
\begin{equation}
\scalebox{0.8}{$\displaystyle
\begin{split}
\int_\Omega\kappa|\nabla(u-u_H)|^2 &\leq \frac{C_1}{\Lambda_*} + C_2H^2 + C 
\sum_i\int_{\omega_i}\kappa|\nabla\chi_i|^2(I^{\omega_i}u-I_0^{\omega_i}u)^2
+C\sum_i\int_{\omega_i}\kappa|\nabla\chi_i|^2(I_0^{\omega_i}u-I_1^{\omega_i}u)^2
\\&  \ \ \ +C \sum_i\int_{\omega_i}\kappa|\nabla\chi_i|^2(I_1^{\omega_i}u-I_2^{\omega_i}u)^2.
\end{split}$}
\end{equation}
We assume that $\mathrm{\Xi}_{train}$ is large 
and define
\[\varepsilon^2_r= \sum_{i,j} \cfrac{1}{|\mathrm{\Xi}_{train}|}\sum_s \|\phi^{\omega_i}_j(\mu_s)-\hat{\phi}^{\omega_i}_j(\mu_s)\|_{L^2}^2 = \sum_{i,j} \mathbb{E}( \|\phi^{\omega_i}_j-\hat{\phi}^{\omega_i}_j\|_{L^2}^2) - \delta. \]
As $|\mathrm{\Xi}_{train}| \rightarrow \infty$, we derive that $\delta \rightarrow 0 $ almost surely provided by the law of large number.
The analysis of POD expansion show that $\varepsilon_{r}^2 = \sum_{k = r+1}^{\infty} \lambda_k$. We denote $S = \sum_{k=1}^{\infty}\lambda_k$ and take $r$ satisfying $\frac{\sum_{k=1}^r\lambda_k}{\sum_{k=1}^{\infty}\lambda_k} \geq 1- \gamma$.
As $r \rightarrow \infty$, we get  $\gamma \rightarrow 0$.
We denote $\| I^{\omega_i}u\|_{s,\text{max}} = \text{max}_{\omega_i}\|I^{\omega_i}u\|_s$ and assume that it is bounded, then we have
\begin{equation}
\scalebox{0.8}{$\displaystyle
\begin{split}
\sum_iE\left( \int_{\omega_i}\kappa|\nabla\chi_i|^2(I^{\omega_i}u-I_0^{\omega_i}u)^2\right) 
&\leq \sum_i\int_P\int_{\omega_i}\kappa H^{-2}(I^{\omega_i}u-I_0^{\omega_i}u)^2dP(\mu)
\\ &=   \sum_i \int_P\int_{\omega_i}\widetilde{\kappa}\Big(\sum_{j=1}^{l_i}c_j(\phi^{\omega_i}_j-\hat{\phi}^{\omega_i}_j)\Big)^2dP(\mu)
    \\ &\leq C \widetilde{\kappa}_{\text{max}}\sum_i\int_P\int_{\omega_i}\Big(\sum_{j=1}^{l_i}c^2_j\Big)\Big(\sum_{j=1}^{l_i}(\phi^{\omega_i}_j-\hat{\phi}^{\omega_i}_j)^2\Big)dP(\mu).
\end{split}$}
\end{equation}
We write $\|I^{\omega_i}u\|_s^2 = \int_{\omega_i}\widetilde{\kappa}\left(\sum_{h=1}^{l_i}c_h\phi^{\omega_i}_h\right)\left(\sum_{j=1}^{l_i}c_j\phi^{\omega_i}_j\right)dx = \sum_{j=1}^{l_i}c_j^2$, then we have
\begin{equation}
\scalebox{0.8}{$\displaystyle
\begin{split}
\sum_iE\left( \int_{\omega_i}\kappa|\nabla\chi_i|^2(I^{\omega_i}u-I_0^{\omega_i}u)^2\right) &\leq C \widetilde{\kappa}_{\text{max}} \sum_{i,j}\int_P\int_{\omega_i}\| I^{\omega_i}u\|_s^2 (\phi^{\omega_i}_j-\hat{\phi}^{\omega_i}_j)^2dP(\mu)
\\ &\leq C \widetilde{\kappa}_{\text{max}}\|I^{\omega_i}u\|_{s,\text{max}}^2 \left(\varepsilon_r^2 + |\delta| \right)
\\ &\leq \hat{C} \widetilde{\kappa}_{\text{max}} (S \gamma + |\delta|)
\\ &\leq C_3 \widetilde{\kappa}_{\text{max}}\gamma + C_3 \widetilde{\kappa}_{\text{max}}|\delta|.
\end{split}$}
\end{equation}
The last step is provided by $S =  \sum_{k=1}^{\infty} \lambda_k \leq C$ for $C$ a constant.
From gPC approximation, we have $\hat{\phi}^{\omega_i}_{j,1} = \sum_{k=1}^{r}\widetilde{a}^{\omega_i,j}_k\psi_k = \sum_{k=1}^{r} \sum c^{\omega_i,j,k}_{\bm{\alpha}} \eta_{\bm\alpha}(\mu) \psi_k$. 
Thus, $\hat{\phi}_j^{\omega_i} - \hat{\phi}_{j,1}^{\omega_i} = \sum_{k=1}^{r}(a_k^{\omega_i,j} - \widetilde{a}_k^{\omega_i,j})\psi_k = \sum_{k=1}^{r}\psi_k\epsilon_{N}^{\omega_i,j,k}$, where     $N$ is the terms of gPC expansion. By integration, we have
\begin{equation}
\scalebox{0.8}{$\displaystyle
\begin{split}
\sum_iE\left( \int_{\omega_i}\kappa|\nabla\chi_i|^2(I_0^{\omega_i}u-I_1^{\omega_i}u)^2\right) &\leq \sum_i\int_P\int_{\omega_i}\kappa H^{-2}(I_0^{\omega_i}u-I_1^{\omega_i}u)^2dP(\mu)
\\ &=   \sum_i \int_P\int_{\omega_i}\widetilde{\kappa}\Big(\sum_{j=1}^{l_i}c_j(\hat{\phi}^{\omega_i}_j-\hat{\phi}^{\omega_i}_{j,1})\Big)^2dP(\mu)
\\ &\leq C  \sum_{i,j}\int_P\int_{\omega_i} \widetilde{\kappa}\|I^{\omega_i}u \|_s^2 (\hat{\phi}^{\omega_i}_j-\hat{\phi}^{\omega_i}_{j,1})^2dP(\mu)
\\ &\leq C \widetilde{\kappa}_{\text{max}}\|I^{\omega_i}u\|_{s,\text{max}}^2 \sum_{i,j}\int_{\omega_i}\int_P \left(\sum_{k=1}^r\psi_k^2\right) \left(\sum_{k=1}^r(\epsilon_N^{\omega_i,j,k})^2\right) dP(\mu)
\\ &\leq  \hat{C} \widetilde{\kappa}_{\text{max}}\sum_{i,j,k}\int_{\omega_i}\psi_k^2 dx\int_P(\epsilon_N^{\omega_i,j,k})^2dP(\mu)
\\ &\leq C_4 \widetilde{\kappa}_{\text{max}}\sum_{i,j,k}\int_P(\epsilon_N^{\omega_i,j,k})^2dP(\mu).
\end{split}$}
\end{equation}
The last ineuality is provided by that $\int_{\omega_i}\psi_k^2 dx$ is bounded and the gPC expansions converge in mean square sense, i.e. $\int_P(\epsilon_N^{\omega_i,j,k})^2dP(\mu) \rightarrow 0$ as N goes infinity. 
Finally, we discuss the error between gPC approximation function and predicted function. In this case, we take least square method and other interpolation methods can be considered. Following with $\widetilde{a}^{\omega_i,j}_k \in \text{span}_{\bm{\alpha}\in \mathbb{N}_0^{\text{M}}}\{\eta_{\bm{\alpha_1}}(\mu),...,\eta_{\bm{\alpha_{\text{M}}}}(\mu) \}$, when the number of samples is greater than or equal to the number of polynomials, we can get the exact solution, that is $E(\hat{\phi}^{\omega_i}_{j,1}-\hat{\phi}^{\omega_i}_{j,2}) = 0$. Consequently, $E(I_1^{\omega_i}u - I_2^{\omega_i}u) = 0$.
Thus we have the following convergence rate for gPC-GMsFEM
\begin{equation}
\scalebox{0.8}{$\displaystyle
\begin{split}
E\left(\int_\Omega\kappa|\nabla(u-u_H)|^2\right) &\leq \frac{C_1}{\Lambda_*} + C_2H^2 +
C_3 \widetilde{\kappa}_{\text{max}}\gamma + C_3 \widetilde{\kappa}_{\text{max}}|\delta| + C_4\widetilde{\kappa}_{\text{max}}\sum_{i,j,k}\int_P(\epsilon_N^{\omega_i,j,k})^2dP(\mu).
\end{split}$}
\end{equation}
\end{proof}

Note that in the analysis above, the convergence rate is proportional to $\frac{1}{\Lambda_*}$ if we take $H \rightarrow 0$ and $r \rightarrow \infty$ as well as $N \rightarrow \infty$ . In practice we cannot guarantee this, but we still have good approximation if we correctly choose the number of eigenfunctions as well as the terms of POD and gPC expansions.

\section{Numerical results}
\label{sec:Numerical}
In this section, we present several numerical examples to illustrate the applicability of the randomized GMsFEM for solving parameterized elliptic partial differential equations. In subsection~\ref{sub:6.1}, we present the normal GMsFEM, gPC-GMsFEM and GPR-GMsFEM and compare them for elliptic PDEs with one dimensional parameter. Moreover, we consider the relationship between the number of basis functions and the convergence results in gPC-GMsFEM. We fix the number of fine and coarse grids, gradually increase the number of basis functions on the coarse grid, and observe the relationship between the dimensionality of the coarse grid and the convergence result. In subsection~\ref{sub:6.3}, we equip $\kappa$ with multiple high contrasts and present the performance of gPC-GMsFEM and GMsFEM.

Before presenting the individual examples, we describe the computational domain, that is, spatial domain $\mathrm{\Omega}=\left[0,1\right]\times[0,1]$, and $\Omega_p$ is a parameter space that satisfies the standard normal Gaussian distribution.
Let $\kappa\left(x;\mu\right):\mathrm{\Omega}\times P\rightarrow R$ be a diffusion coefficient function. We assume that $\mu$ is a parameter sampled from a Gaussian distribution defined over the parameter space 
$\Omega_p$. We consider the following model elliptic equation for numerical computation,
\begin{equation}
    -\mathrm{\nabla}\cdot\left(\kappa\left(x;\mu\right)\mathrm{\nabla u}\right)=f \ \text{in} \ \mathrm{\Omega},\ u=p\ \text{on}\ \partial\Omega.
\end{equation}

\subsection{Numerical results for one-dimensional periodic parameter}
\label{sub:6.1}
\begin{figure}[htbp]
    \centering
    \includegraphics[width=0.9\textwidth, height=0.35\textheight]{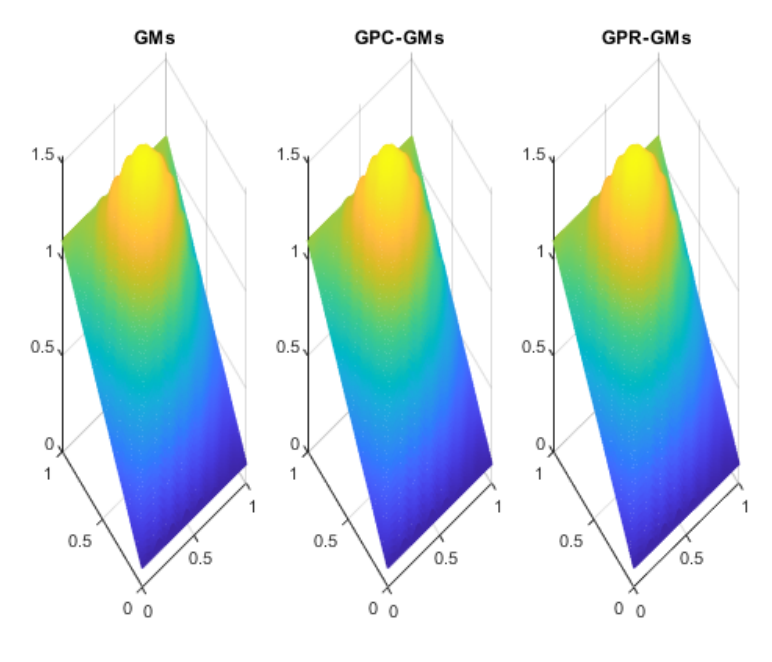}
    \caption{Numerical solutions of GMsFEM (left), gPC-GMsFEM (middle), and GPR-GMsFEM (right) for the permeability field depicted in (\ref{kp}). Here $n_x=100$, $N_x=5$.}
    \label{fig:1}
\end{figure}

In this subsection, we will consider one numerical example. 
We assume that
\begin{subequations}
    \begin{gather}
        \kappa\left(x;\mu\right)=\left(x^2y+\left(3+2.8\sin{15\pi\left(x-y\right)}\right)^{-1}\right)\left(\mu+\mu^2\right),\label{kp}\\
        f=\pi^2\left(2+y\right)\sin{\pi x}\sin{\pi y}+2\pi^2x\cos{\pi x}\cos{\pi y},\\
        p=\sin{\pi x}\sin{\pi y}+y+0.1,
    \end{gather}
\end{subequations}
where $\mu \in \mathbb{R}$. In the offline stage, we use the space composed of fine grid basis functions as a snapshot space. When solving the generalized eigenvalue problem, we take the first $l_i$ smallest eigenvectors and multiply them by standard multiscale basis functions to generate a coarse space. The sample points were randomly sampled and the gPC prediction functions were obtained using the least squares method. We take $\mu+\mu^2>0$ to ensure the positive definiteness of $\kappa\left(x;\mu\right)$. In the online stage, the feature vector is predicted based on the new parameter $\mu_\ast=0.6$, and multiplied by standard multiscale basis functions to generate the online coarse space, obtaining the projection matrix. At the same time, the right-hand space is generated based on the values of $p$ and $f$.  
We take the fine grid as $n_x\times n_y=100\times 100$ and the coarse grid as $N_x\times\ N_y=5\times 5$. The obtained image is shown in Figure~\ref{fig:1}.

We present the respective errors between GMsFEM, gPC-GMsFEM, and GPR-GMsFEM compared to the fine scale Finite Element Method (FEM). We define 
the energy norm as $||u||_E=\sqrt{\int_{\mathrm{\Omega}}{\kappa\nabla u\cdot\mathrm{\nabla u}}}$. The $L^2$ norm is $||u||_{L^2}=\sqrt{\int_{\mathrm{\Omega}}\kappa\left|u\right|^2}$. The errors between GMsFEM, gPC-GMsFEM, and GPR-GMsFEM compared to the fine scale FEM are presented in the Table~\ref{tab:1}, where the N-basis denotes the number of basis functions we take for each $\omega_i$.

\begin{table}[htbp]
{\footnotesize
  \caption{Relative $L^2$ and energy errors for the permeability field depicted in (\ref{kp}). Here $n_x=100$, $N_x=5$, N-basis = 5.}  \label{tab:1}
\begin{center}
  \begin{tabular}{|c|c|c|c|} \hline
    Error & GMsFEM & gPC-GMsFEM & GPR-GMsFEM \\ \hline
    Relative $L^2$-norm error & 0.0058 & 0.0058 & 0.0058 \\
    Relative energy norm error & 0.0967 & 0.0967 & 0.0967\\ \hline
  \end{tabular}
\end{center}
}
\end{table}

In Figure~\ref{fig:1}, we present the solutions obtained using GMsFEM, gPC-GMsFEM, and GPR-GMsFEM. All three methods provide a good approximation for the original problem, thereby demonstrating the feasibility of the method. In Table~\ref{tab:1}, we report the relative errors of different types. The gPC-GMsFEM and GPR-GMsFEM both yield a good approximation for the original problem, leading to comparable convergence results. This can be interpreted as the error of GMsFEM itself being larger than the error of the predictors, so the error in this problem is dominated by the error in GMsFEM. Next we take the fine grid as $n_x\times n_y=100\times100$ and the coarse grid as $N_x\times N_y=10\times10$. The relative errors obtained are presented in Table~\ref{tab:2}.

\begin{table}[htbp]
{\footnotesize
  \caption{Relative $L^2$ and energy errors for the permeability field depicted in (\ref{kp}). Here $n_x=100$, $N_x=10$, N-basis = 5.}  \label{tab:2}
\begin{center}
  \begin{tabular}{|c|c|c|c|} \hline
    Error & GMsFEM & gPC-GMsFEM & GPR-GMsFEM \\ \hline
    Relative $L^2$-norm error & 0.0039 & 0.0039 & 0.0039 \\
    Relative energy norm error & 0.0922 & 0.0922 & 0.0923\\ \hline
  \end{tabular}
\end{center}
}
\end{table}

As shown in Table~\ref{tab:2}, with finer grids, the relative errors of randomized GMsFEM are narrowing, which is consistent with the actual operation. It can also be seen that the method of using predictors have a comparable approximation due to the same reason. At the same time, we found in the experiment that the performance of GPR-GMsFEM deteriorates significantly when the training data contains outliers or insufficient samples. Poor sampling data can lead to significant discrepancies between the predicted and actual values of the GPR-GMsFEM.

For gPC-GMsFEM, we fix the fine and coarse grid sizes. By gradually increasing the number of basis functions on the coarse grid, we analyze how the coarse grid dimensionality affects the convergence results. As the dimensionality of the coarse space increases, the relative errors of gPC-GMsFEM solutions decrease and the convergence effect gradually improves. The relationship between the dimensionality of a single coarse grid and the convergence result is shown in the Table~\ref{tab:3}, where the content in brackets indicates the total degrees of freedom.

\begin{table}[htbp]
{\footnotesize
  \caption{Convergence history for the permeability field depicted in (\ref{kp}) and for the case with one initial basis applying gPC-GMsFEM. Here $n_x=100$, $N_x=5$.} \label{tab:3}
\begin{center}
  \begin{tabular}{|c|c|c|c|c|c|} \hline
    N-basis & $L^2$Error & Energy error & N-basis & $L^2$Error & Energy error\\ \hline
    1(36)  & 0.0332 & 0.2654 & 6(216)  & 0.0041 & 0.0844\\ 
    2(72)  & 0.0182 & 0.1832 & 7(252)  & 0.0032 & 0.0736\\ 
    3(108)  & 0.0158 & 0.1688 & 8(288)  & 0.0027 & 0.0678\\ 
    4(144)  & 0.0085 & 0.1192 & 9(324)  & 0.0025 & 0.0640\\ 
    5(180)  & 0.0058 & 0.0967 & 10(360)  & 0.0022 & 0.0599\\\hline 

  \end{tabular}
\end{center}
}
\end{table}

\subsection{Numerical results for $\kappa$ with high contrast}
\label{sub:6.3}
In this subsection, we will consider $\kappa(x;\mu)$ with high contrasts. We present the field of the contrast $\tau=10^4$, where $\tau$ represent the ratio between the highest region and the lowest region. We still assume that $\mu \in \mathbb{R}^{\text{M}}$, satisfying the normal Gaussian distribution and we inherit the assumptions as subsection~\ref{sub:6.1} did except for $\kappa$. We take the fine grid as $n_x\times n_y=100\times 100$ and the coarse grid as $N_x\times\ N_y=5\times 5$. 

\subsubsection{Case 1: Equation with one high contrast}\label{sec:case1}
\begin{figure}[htbp]
    \centering
    \vspace{-90pt}
    \includegraphics[width=0.70\textwidth]{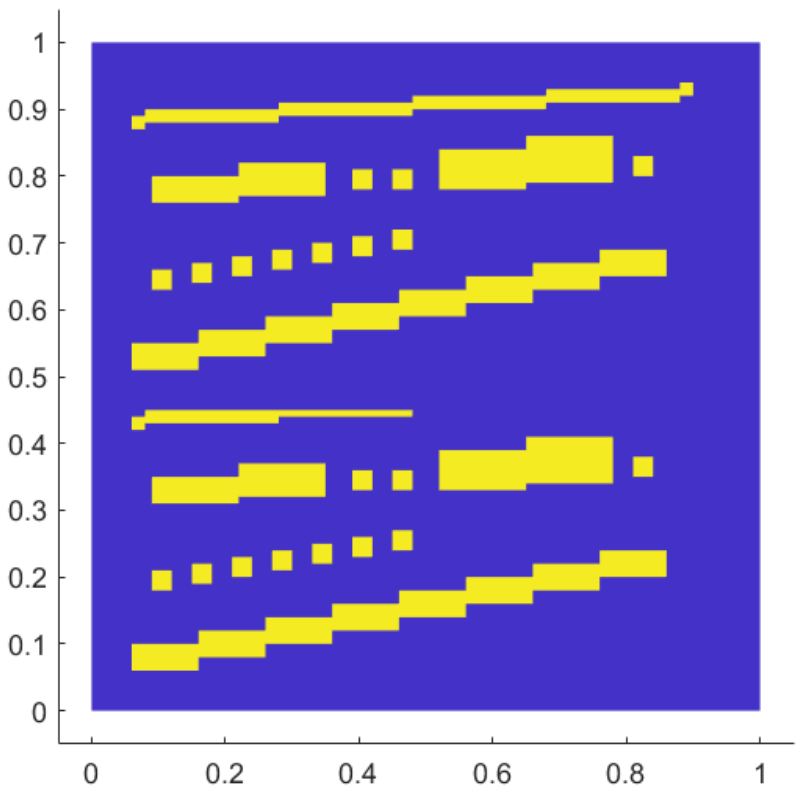}
    \vspace{-80pt}
    \caption{Permeability field $\kappa$, $\tau=10^4$.}
    \label{fig:2}
\end{figure}

\begin{figure}[htbp]
    \centering
    \includegraphics[width=0.9\textwidth, height=0.35\textheight]{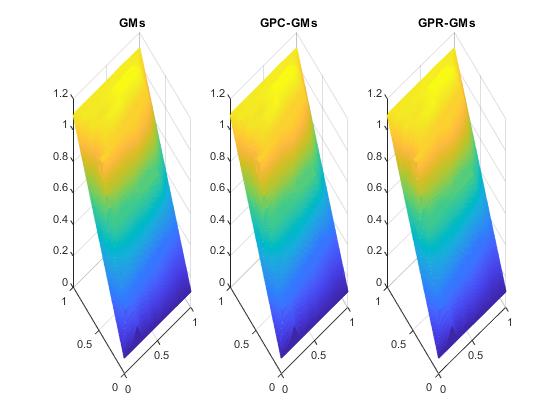}
    \caption{Numerical solutions of GMsFEM (left), gPC-GMsFEM (middle), and GPR-GMsFEM (right) for the permeability field depicted in Figure~\ref{fig:2}.}
    \label{fig:3}
\end{figure}

we begin with a high contrast case, as shown in Figure~\ref{fig:2}, where \( Q = 1 \) and the permeability field is given by \( \kappa(x;\mu) = \Theta_1(\mu)\kappa_1(x) \).
In Figure~\ref{fig:3}, we compare the numerical solutions obtained using GMsFEM (left), gPC-GMsFEM (middle), and GPR-GMsFEM (right) through numerical visualizations.

Correspondingly, in Table~\ref{tab:4} we show the relative errors of the three methods. We calculate the errors especially the relative $L^2$ error several times, and find that in most cases the gPC-GMsFEM performs slightly better than GMsFEM. The performance of GPR-GMsFEM seems not good compared with gPC-GMsFEM and it costs a long time. The average time for solving this case using gPC-GMsFEM as well as GMsFEM is presented in Table~\ref{tab:111} and we can see that gPC-GMsFEM has a better computational efficiency.  To show the numerical convergence of gPC-GMsFEM, we present Table~\ref{tab:5} of the relationship between the dimensionality of a single coarse grid and the convergence result. After testing, we found that as the number of local basis functions increases, the convergence results of GMsFEM and gPC-GMsFEM are basically consistent, only a slight fluctuation of gPC-GMsFEM has and can be ignored.

\begin{table}[htbp]
{\footnotesize
  \caption{Relative $L^2$ and energy errors for the permeability field depicted in Figure~\ref{fig:2}. Here N-basis = 5.}  \label{tab:4}
\begin{center}
  \begin{tabular}{|c|c|c|c|} \hline
    Error & GMsFEM & gPC-GMsFEM & GPR-GMsFEM \\ \hline
    Relative $L^2$-norm error & 0.0412 & 0.0411 & 0.0412 \\
    Relative energy norm error & 0.2985 & 0.2984 & 0.2987\\ \hline
  \end{tabular}
\end{center}
}
\end{table}

\begin{table}[htbp]
{\footnotesize
    \caption{Average CPU time of solving the permeability field depicted in Figure~\ref{fig:2} in the online stage.} \label{tab:111}
\begin{center}
    \begin{tabular}{|c|c|c|c|c|c|} \hline
        Algorithm  & N-basis & CPU Time & Algorithm & N-basis &CPU Time  \\ \hline
                   & 1(36)   & 0.080s  &         & 1(36)   & 0.626s \\
                   & 3(108)  & 0.096s  &         & 3(108)  & 0.652s \\
                   & 5(180)  & 0.115s  &         & 5(180)  & 0.714s \\
    gPC-GMsFEM     & 7(252)  & 0.133s  &  GMsFEM & 7(252)  & 0.882s \\
                   & 9(324)  & 0.147s  &         & 9(324)  & 0.925s \\
                   & 11(396)  & 0.175s &         & 11(324)  & 0.992s\\
                   & 13(468)  & 0.209s &         & 13(324)  & 1.060s\\ \hline
    \end{tabular}
\end{center}
}
\end{table}

\begin{table}[htbp]
{\footnotesize
    \caption{Convergence history for the permeability field depicted in Figure~\ref{fig:2} and for the case with one initial basis applying gPC-GMsFEM.}
  \label{tab:5}
\begin{center}
  \begin{tabular}{|c|c|c|c|c|c|} \hline
    N-basis  & $L^2$Error & Energy error & N-basis  & $L^2$Error & Energy error  \\ \hline 
    1(36)  & 0.3917 & 2.1524 & 8(288)  & 0.0066 & 0.1293\\
    2(72)  & 0.2991 & 0.7586 & 9(324)  & 0.0039 & 0.1101 \\ 
    3(108)  & 0.1463 & 0.6752 & 10(360)  & 0.0036 & 0.1080\\
    4(144)  & 0.1199 & 0.4662 & 11(396)  & 0.0027 & 0.0988 \\ 
    5(180)  & 0.0411 & 0.2984 & 12(432)   & 0.0024 & 0.0937\\  
    6(216)  & 0.0386 & 0.2226 & 13(468)  & 0.0021 & 0.0876\\ 
    7(252)  & 0.0154 & 0.1535 & 14(504)  & 0.0018 & 0.0805 \\\hline  
  \end{tabular}
\end{center}
}
\end{table}

\begin{figure}[htbp]
    \centering
    \vspace{-80pt}
    \includegraphics[width=0.60\textwidth]{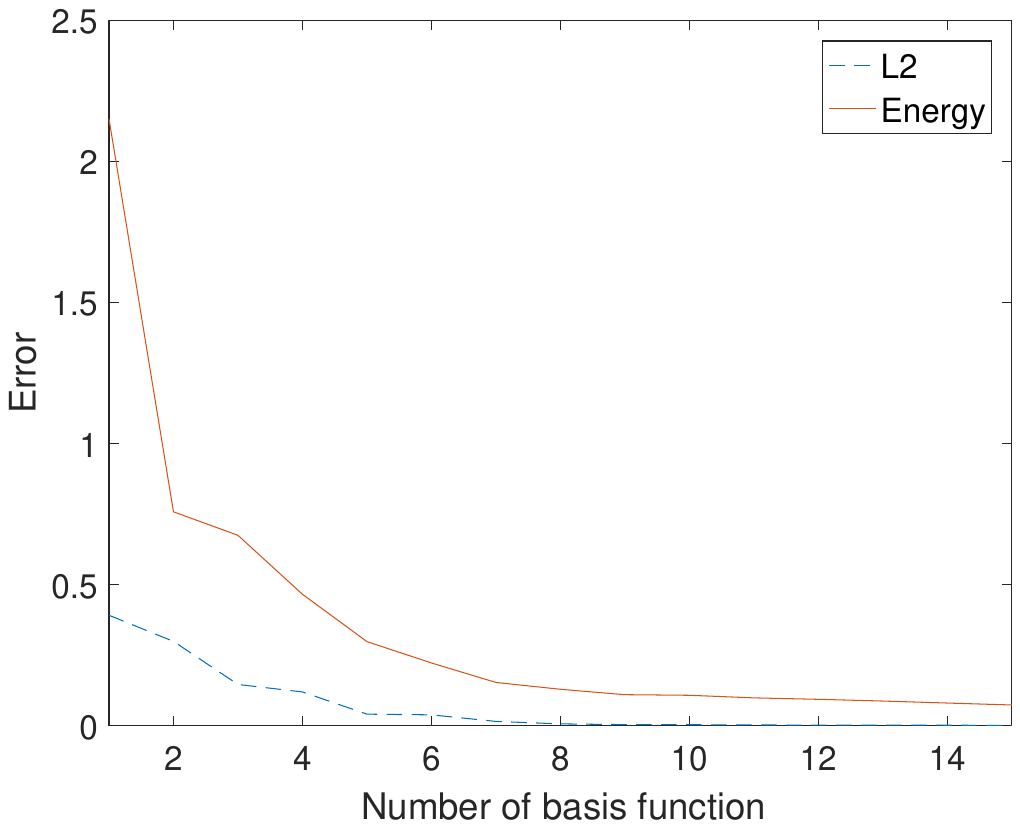}
    \vspace{-80pt}
    \caption{Relative $L^2$ and energy errors with a varying number of basis functions applying gPC-GMsFEM for the permeability field depicted in Figure~\ref{fig:2}.}
    \label{fig:4}
\end{figure}

To further analyze the convergence behavior, Figure~\ref{fig:4} illustrates the relationship between the number of basis functions and the approximation errors using gPC-GMsFEM. The errors decrease approximately exponentially  with the increase of the basis functions. That is because as the channels and inclusions containing high contrast are captured by the former basis functions which have the small eigenvalues, the role of the latter basis functions will no longer be significant.

\subsubsection{Case 2: Equation with two high contrasts}\label{sec:case2} 
\begin{figure}[htbp]
	\centering
	\begin{minipage}{0.48\linewidth}
		\centering
		\includegraphics[width=0.95\linewidth]{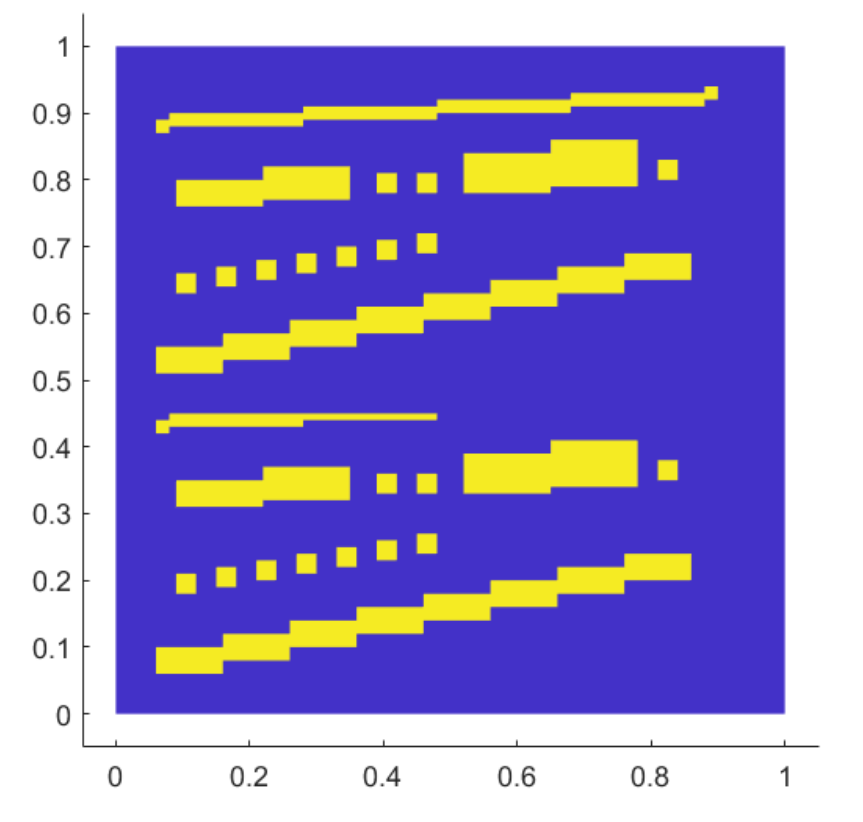}
	\end{minipage}
	\begin{minipage}{0.48\linewidth}
		\centering
		\includegraphics[width=0.95\linewidth]{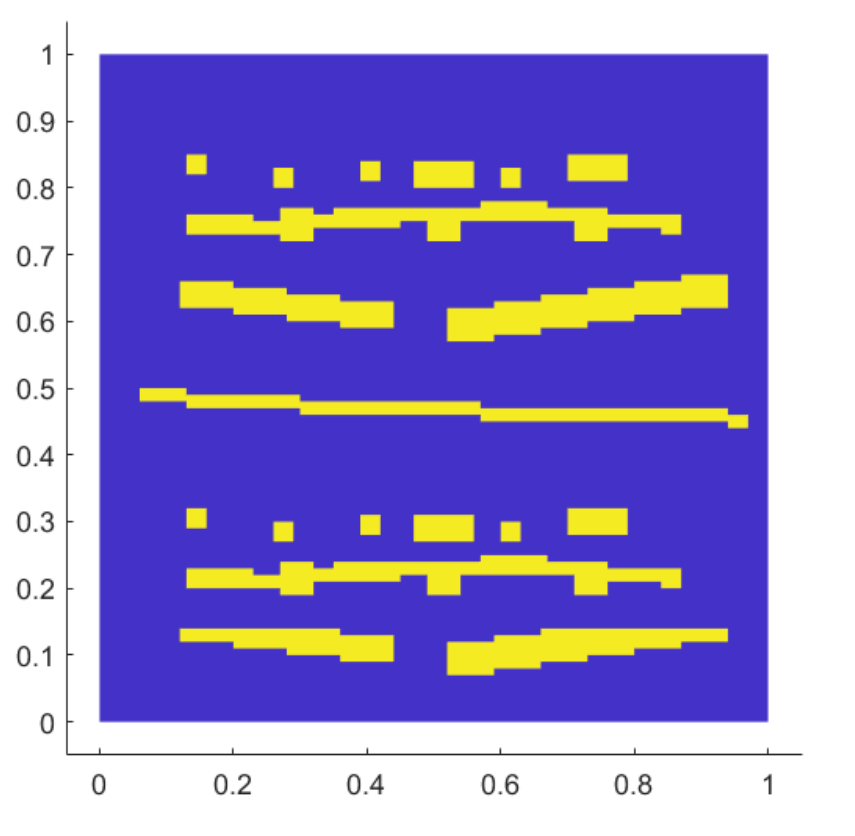}
            \vspace{2pt}
	\end{minipage}
    \caption{Permeability fields. Left: $\kappa_1$ and right: $\kappa_2$, $\tau=10^4$.}
    \label{fig:5}
\end{figure}

\begin{figure}[htbp]
	\begin{minipage}{0.45\linewidth}
		\centerline{\includegraphics[width=\textwidth]{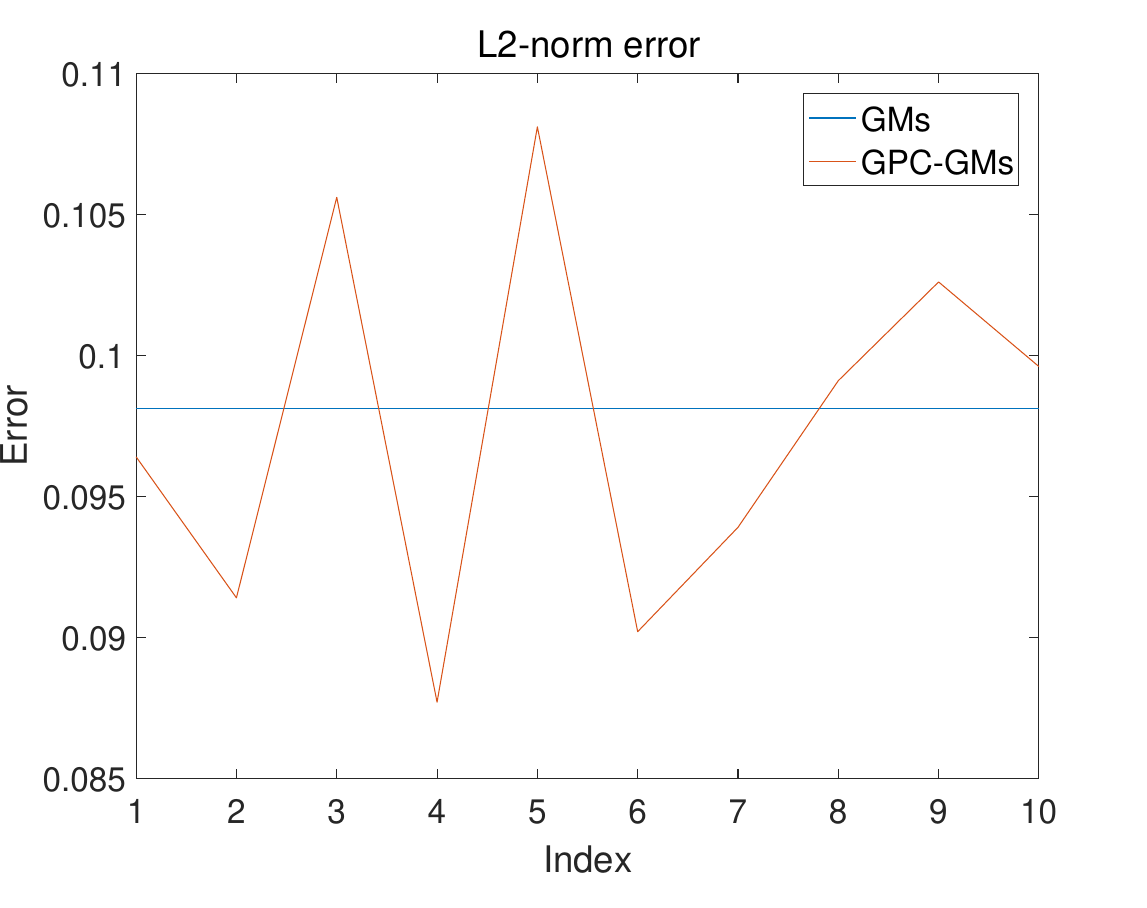}}
	\end{minipage}
	\begin{minipage}{0.45\linewidth}
		\centerline{\includegraphics[width=\textwidth]{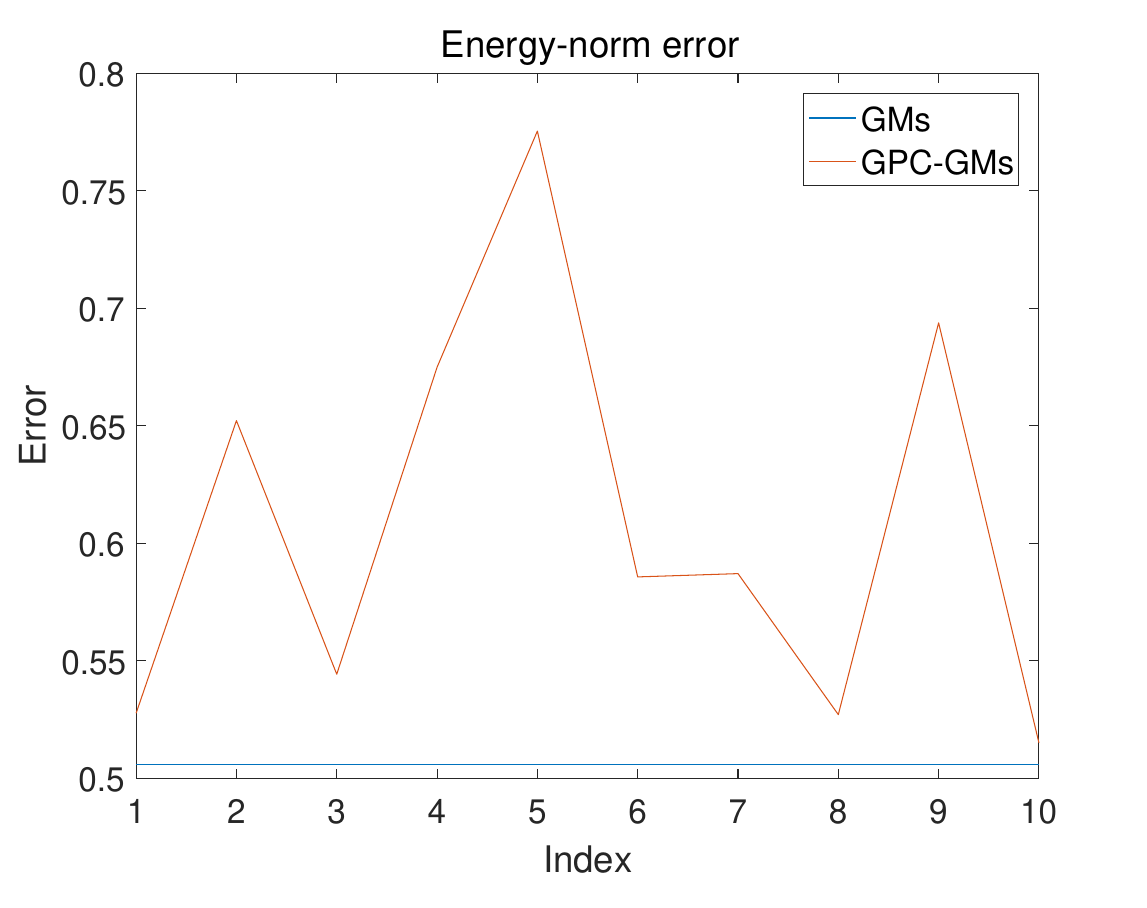}}
	\end{minipage}
	\caption{Relative $L^2$ and energy errors for the permeability fields depicted in Figure~\ref{fig:5} with different predictors. Index indicates that the predictors are trained with different samples. Left: $L^2$ error, and right: energy error.}
	\label{fig:6}
\end{figure}

\begin{figure}[htbp]
    \centering
    \vspace{-80pt}
    \includegraphics[width=0.6\textwidth]{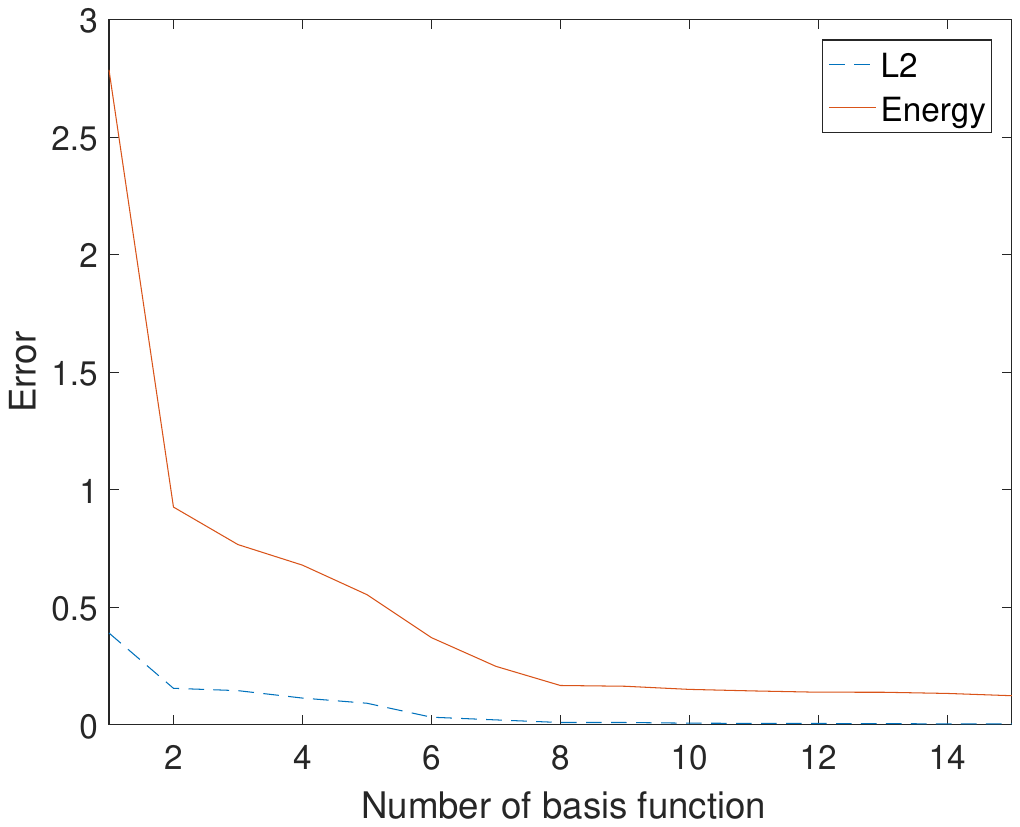}
    \vspace{-80pt}
    \caption{Relative $L^2$ and energy errors with a varying number of basis functions applying gPC-GMsFEM for the permeability fields depicted in Figure~\ref{fig:5}.}
    \label{fig:2k}
\end{figure}

We consider a more complicated situation. Assuming that $\kappa$ to be $\kappa\left(x;\mu\right)=\sum_{q=1}^{Q}{\Theta}_q\left(\mu\right)\kappa_q\left(x\right)$. To illustrate clearly, we take $Q=2$, $\mu=(\mu^1, \mu^2, \dots , \mu^M)\in \mathbb{R}^{\text{M}}$ and obtain $\kappa_1$ as well as $\kappa_2$ given in Figure~\ref{fig:5}. 
The gPC expansion of finitely normal basic random variables should be $\phi\left(\mu\right)=\sum_{i=0}^{N}{c_i\bm\eta_i}\left({\mu}\right)$, where $\bm\eta_i(\mu)=\prod_{j=1}^d{\eta_{\alpha_j}^i(\mu^j)}, \bm\alpha=(\alpha_1,\alpha_2, \dots , \alpha_d)\in \mathbb{N}_0^{\text{M}}$. Due to the complicacy of the parameter coefficient, it will bring about an approximated difficulty in our numerical method. It is reflected in the choice of random sampling, and we note that influenced by the random sample, the approximation of the gPC expansion shows an oscillatory trend as the number of experiments increases. We present a group of plots to clarify this situation in Figure~\ref{fig:6}. We use the error of GMsFEM as a reference, finding that the $L^2$ error of gPC-GMsFEM is comparable to the GMsFEM. Besides, the difference in energy error between gPC-GMsFEM and GMsFEM is not significant. In this case, the GPR-GMsFEM will incur huge computational costs and its representation is worse than the gPC-GMsFEM compared in the same condition, which is shown in Table~\ref{tab:6}, and we take the average of the errors of gPC-GMsFEM of the figure above. In Table~\ref{tab:2k}, we present the relationship between errors and the dimensionality of the coarse grid, that is, the number of basis functions. We found that as the number of basis functions increases, the errors decays exponentially roughly and we present a line plot in Figure~\ref{fig:2k} to show it clearly.

\begin{table}[htbp]
{\footnotesize
  \caption{Relative $L^2$ and energy errors for the permeability fields depicted in Figure~\ref{fig:5}. Here N-basis = 5.}  \label{tab:6}
\begin{center}
  \begin{tabular}{|c|c|c|c|} \hline
    Error & GMsFEM & gPC-GMsFEM & GPR-GMsFEM \\ \hline
    Relative $L^2$-norm error & 0.0981 & 0.0975 & 0.1482 \\
    Relative energy norm error & 0.5059 & 0.6083 & 0.9078\\ \hline
  \end{tabular}
\end{center}
}
\end{table}

\begin{table}[htbp]
{\footnotesize
  \caption{Convergence history for the permeability fields depicted in Figure~\ref{fig:5} and for the case with one initial basis applying gPC-GMsFEM.}  \label{tab:2k}
\begin{center}
  \begin{tabular}{|c|c|c|c|c|c|} \hline
    N-basis  & $L^2$Error & Energy error & N-basis  & $L^2$Error & Energy error \\ \hline 
    1(36)    & 0.3903 & 2.7865 & 8(288)   & 0.0095 & 0.1670 \\
    2(72)    & 0.1550 & 0.9262 & 9(324)   & 0.0095 & 0.1638 \\ 
    3(108)   & 0.1451 & 0.7666 & 10(360)  & 0.0067 & 0.1504\\
    4(144)   & 0.1135 & 0.6791 & 11(396)  & 0.0058 & 0.1438\\ 
    5(180)   & 0.0915 & 0.5536 & 12(432)  & 0.0053 & 0.1384\\  
    6(216)   & 0.0321 & 0.3707 & 13(468)  & 0.0045 & 0.1383\\ 
    7(252)   & 0.0206 & 0.2485 & 14(504)  & 0.0038 & 0.1332\\\hline 

  \end{tabular}
\end{center}
}
\end{table}

\subsubsection{Case 3: Equation with four high contrasts}\label{sec:case3} 
\begin{figure}[htbp]
 \begin{minipage}{0.46\linewidth}
    \vspace{2pt}
    \centerline{\includegraphics[width=\textwidth]{fig/kappa1_zy_1.pdf}}
     \vspace{9pt}
     \centerline{\includegraphics[width=\textwidth]{fig/kappa2_zy2.pdf}}
 \end{minipage}
    \begin{minipage}{0.46\linewidth}
     \centerline{\includegraphics[width=\textwidth]{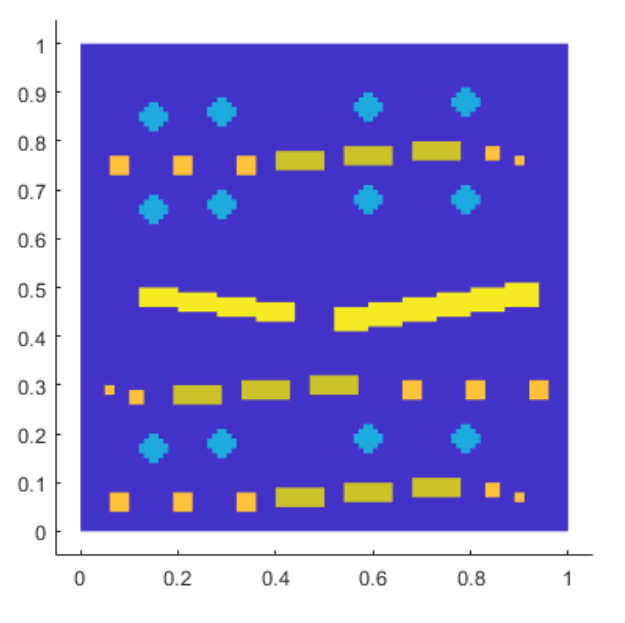}}
     \vspace{3pt}
     \centerline{\includegraphics[width=\textwidth]{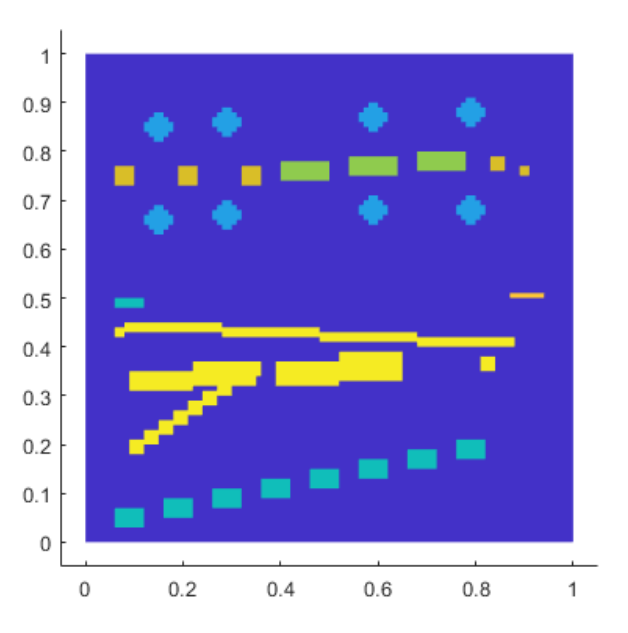}}
 \end{minipage}
\caption{Permeability fields. left: $\kappa_1$, top right: $\kappa_2$, lower left: $\kappa_3$ and lower right: $\kappa_4$, $\tau = 10^4$.}\label{fig:7}
\end{figure}

\begin{figure}[htbp]
	\begin{minipage}{0.475\linewidth}
		\centerline{\includegraphics[width=\textwidth]{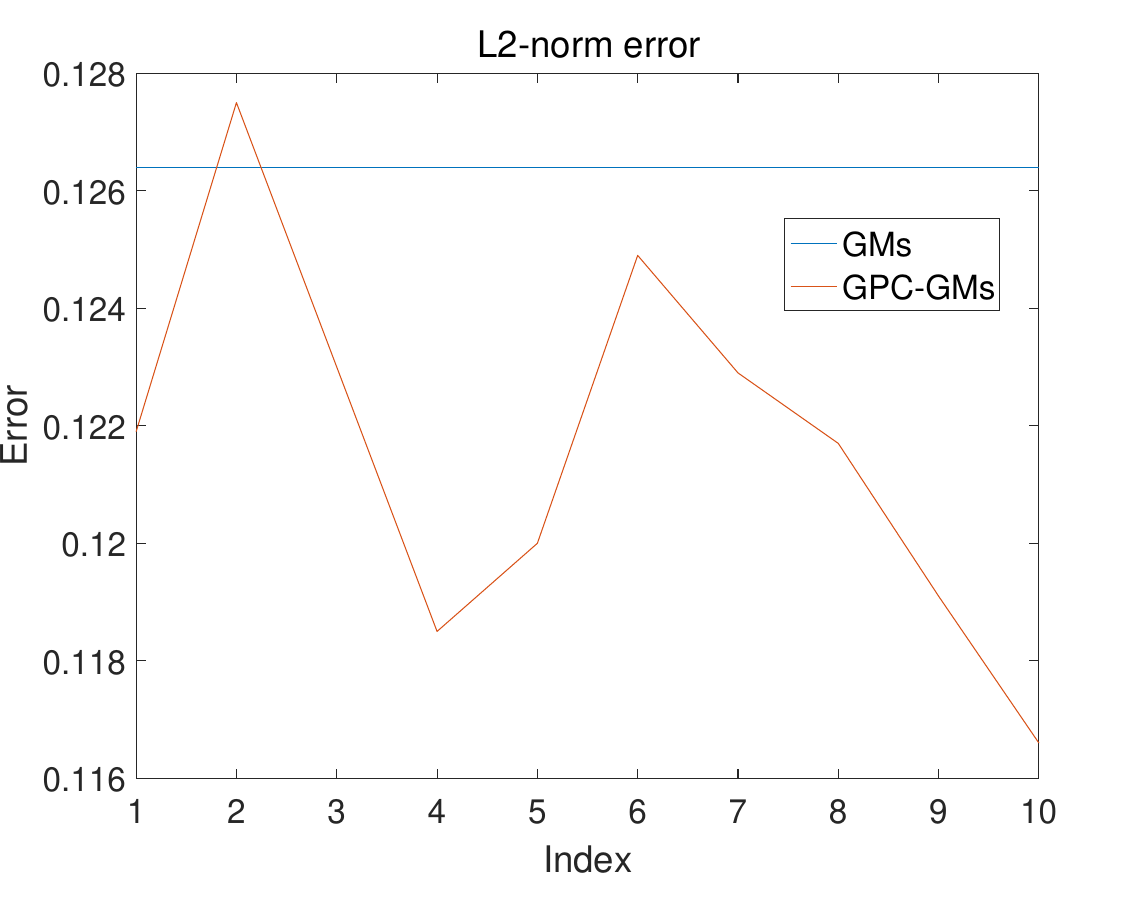}}
	\end{minipage}
	\begin{minipage}{0.475\linewidth}
		\centerline{\includegraphics[width=\textwidth]{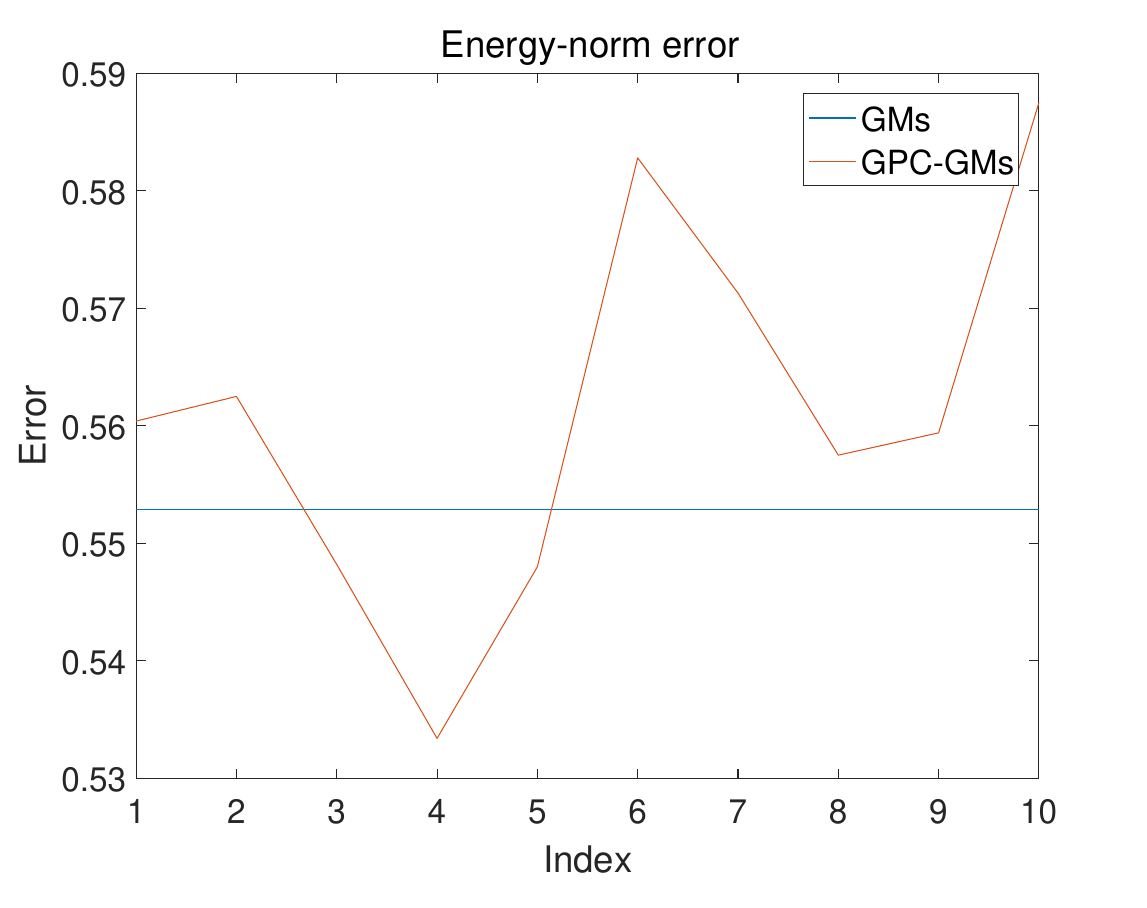}}
	\end{minipage}
	\caption{Relative $L^2$ and energy errors for the permeability fields depicted in Figure~\ref{fig:7} with different predictors. Index indicates that the predictors are trained with different samples. Left: $L^2$ error, and right: energy error.}
	\label{fig:8}
\end{figure}

\begin{figure}[htbp]
    \centering
    \vspace{-80pt}
    \includegraphics[width=0.60\textwidth]{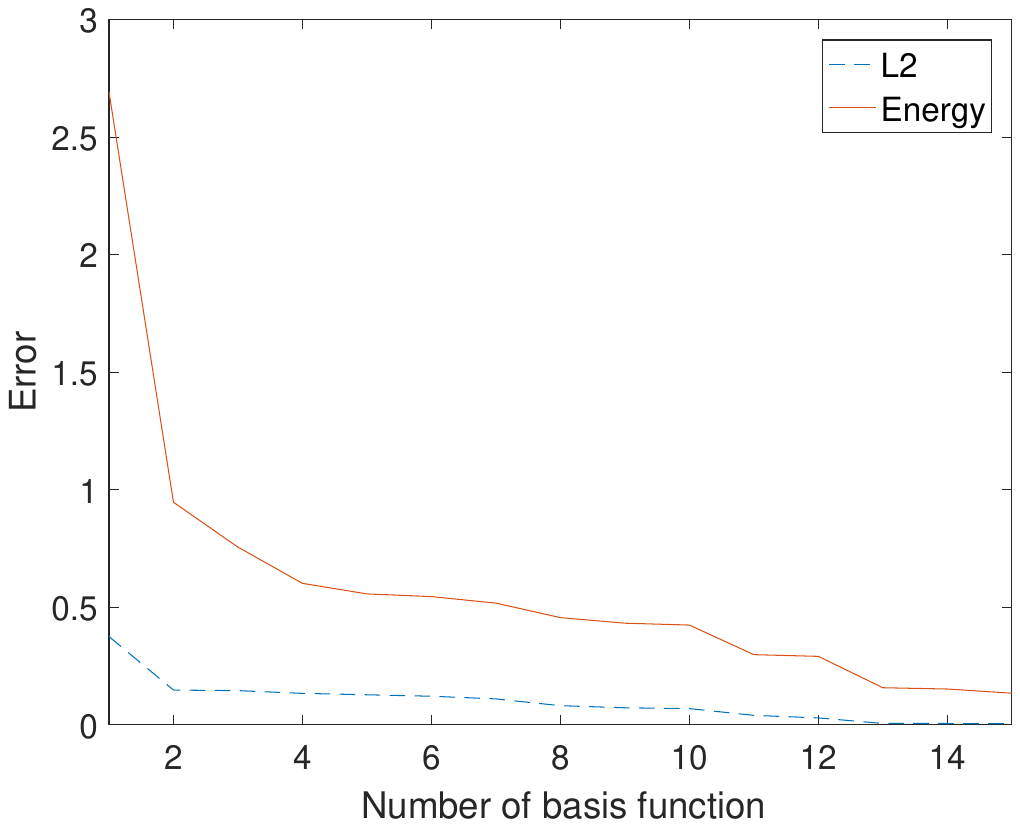}
    \vspace{-80pt}
    \caption{Relative $L^2$ and energy errors with a varying number of basis functions applying gPC-GMsFEM for the permeability fields depicted in Figure~\ref{fig:7}.}
    \label{fig:4k}
\end{figure}

Finally, we take four high contrasts regions and test the performance, that is $Q=4$. We present the images of the four high contrasts coefficients in Figure~\ref{fig:7}. In this situation, it will be a great challenge to approximate the solution and also take huge computational costs. We still take $n_x=100$ and $N_x=5$ and use the previous boundary condition and source function. We present the numerical result in Figure~\ref{fig:8} and we can see that the gPC-GMsFEM method has a great approximation property. Moreover, in the online stage, gPC-GMsFEM will cost less time. In Table~\ref{tab:7}, we present the errors of three different methods. gPC-GMsFEM performs better than GPR-GMsFEM. In order to clarify the efficiency and accuracy of the online stage, we use 50 samples to generate the predictor and use it to construct online basis functions. In Table~\ref{tab:9}, we present the error of different random parameters using five basis functions for each $\omega_i$. Using the predictor, we have successfully reduced computational costs and achieved a good errors approximation. We found that neither GMsFEM nor gPC-GMsFEM approximates the energy error very well. This is because as the complexity of the region increases, multiscale basis functions cannot cover areas with high contrasts well, but this can be effectively improved by adding local multiscale basis functions. In Table~\ref{tab:10}, we show the relationship between error and the choice of basis functions per coarse grid. We found that as the number of basis functions increases, the error decays exponentially as in the case of two high contrasts above. We show our convergence result in Figure~\ref{fig:4k} as before. Finally, in Table~\ref{tab:22}, we present the average time required to solve this case. As the high-contrast part increases, the situation becomes more complicated, which brings about a greater amount of calculation, but our method still has a faster calculation speed.

\begin{table}[htbp]
{\footnotesize
  \caption{Relative $L^2$ and energy errors for the permeability fields depicted in Figure~\ref{fig:7}. Here N-basis = 5.}  \label{tab:7}
\begin{center}
  \begin{tabular}{|c|c|c|c|} \hline
    Error & GMsFEM & gPC-GMsFEM & GPR-GMsFEM \\ \hline
    Relative $L^2$-norm error & 0.1264 & 0.1216 & 0.1787 \\
    Relative energy norm error & 0.5529 & 0.5611 & 0.9100\\ \hline
  \end{tabular}
\end{center}
}
\end{table}

\begin{table}[htbp]
{\footnotesize
  \caption{Relative $L^2$ and energy errors for the permeability fields depicted in Figure~\ref{fig:7} with a varying number of parameters. Here N-basis = 5.}  \label{tab:9}
\begin{center}
  \begin{tabular}{|c|c|c|c|c|} \hline
    $\mu$  & $L^2$ (gPC) & $L^2$ (GMs) & Energy  (gPC) & Energy  (GMs) \\ \hline 
    (0.1,0.1,0.1,0.1)  & 0.2269 & 0.1792 & 0.6817 & 0.5251 \\ 
    (0.2,0.2,0.2,0.2)  & 0.2020 & 0.1572 & 0.6660 & 0.5680\\ 
    (0.3,0.3,0.3,0.3)  & 0.1455 & 0.1889 & 0.8182 & 0.5793\\ 
    (0.4,0.4,0.4,0.4)  & 0.1369 & 0.1383 & 0.5897 & 0.5764\\ 
    (0.5,0.5,0.5,0.5)  & 0.1320 & 0.1333 & 0.5844 & 0.5689\\ 
    (0.6,0.6,0.6,0.6)  & 0.1284 & 0.1298 & 0.5777 & 0.5604\\ 
    (0.7,0.7,0.7,0.7)  & 0.1257 & 0.1271 & 0.5709 & 0.5524\\  
    (0.8,0.8,0.8,0.8)  & 0.1243 & 0.1250 & 0.5994 & 0.5451\\ 
    (0.9,0.9,0.9,0.9)  & 0.1228 & 0.1233 & 0.5951 & 0.5388\\ 
    (0.5,0.4,0.3,0.2)  & 0.1389 & 0.1406 & 0.6118 & 0.5782\\ 
    (0.4,0.6,0.5,0.1)  & 0.1344 & 0.1338 & 0.6212 & 0.5772\\ 
    (0.6,0.7,0.5,0.9)  & 0.1265 & 0.1267 & 0.6038 & 0.5499\\ \hline  
  \end{tabular}
\end{center}
}
\end{table}

\begin{table}[htbp]
{\footnotesize
  \caption{Convergence history for the permeability fields depicted in Figure~\ref{fig:7} and for the case with one initial basis applying gPC-GMsFEM.} \label{tab:10}
\begin{center}
  \begin{tabular}{|c|c|c|c|c|c|} \hline
    N-basis & $L^2$Error & Energy error & N-basis & $L^2$Error & Energy error\\ \hline
    1(36)   & 0.3752 & 2.6931 & 8(288)  & 0.0811 & 0.4558\\ 
    2(72)   & 0.1476 & 0.9467 & 9(324)  & 0.0722 & 0.4321\\ 
    3(108)  & 0.1451 & 0.7557 & 10(360) & 0.0686 & 0.4241\\
    4(144)  & 0.1335 & 0.6017 & 11(396) & 0.0400 & 0.2982\\
    5(180)  & 0.1275 & 0.5567 & 12(432) & 0.0291 & 0.2910\\
    6(216)  & 0.1213 & 0.5453 & 13(468) & 0.0054 & 0.1575\\
    7(252)  & 0.1099 & 0.5174 &14(504) & 0.0052 & 0.1520\\\hline 

  \end{tabular}
\end{center}
}
\end{table}

\begin{table}[htbp]
{\footnotesize
    \caption{Average CPU time of solving the permeability field depicted in Figure~\ref{fig:7} in the online stage.} \label{tab:22}
\begin{center}
    \begin{tabular}{|c|c|c|c|c|c|} \hline
        Algorithm  & N-basis & CPU Time & Algorithm & N-basis &CPU Time  \\ \hline
                   & 1(36)   & 0.084s  &         & 1(36)   & 0.850s \\
                   & 3(108)  & 0.110s  &         & 3(108)  & 0.876s \\
                   & 5(180)  & 0.142s  &         & 5(180)  & 0.912s \\
    gPC-GMsFEM     & 7(252)  & 0.183s  &  GMsFEM & 7(252)  & 0.962s \\
                   & 9(324)  & 0.221s  &         & 9(324)  & 1.012s \\
                   & 11(396)  & 0.261s &         & 11(324)  & 1.071s\\
                   & 13(468)  & 0.303s &         & 13(324)  & 1.217s\\ \hline
    \end{tabular}
    \label{tab:time}
\end{center}
}
\end{table}

\section{Conclusions}
\label{sec:conclusions}

In this paper, we study a randomized GMsFEM using predictors designed to construct online basis functions. The main idea of the proposed method is to construct multiscale online basis functions given a new parameter by multiplying standard multiscale basis functions with predicted eigenfunctions provided by predictors in each coarse domain. We constructed a reduced matrix via SVD in the snapshot space and use it as the training data to generate different predictors. The overall procedure is computed in the offline stage and predictors are stored as matrices. During the online stage, for any new parameter, we avoid recomputing the spectral decomposition and generate the online space via different predictors to solve the equation online. We tested our approach on several examples, particularly multiscale problems involving multiple high contrasts permeability fields with random coefficients. Our numerical results show that gPC-GMsFEM has a good approximation for above problems which pose challenges for other methods. Our convergence study shows that the convergence behavior of gPC-GMsFEM is proportional to $\frac{1}{\Lambda_*}$, where $\Lambda_*$ is the smallest eigenvalue not included in the basis space for all $\mu$.

\appendix
\section{The predictor via GPR} 
\label{GPR}
In this section, we present an approach that constructs the predictor via Gaussian Process Regression (GPR). A Gaussian Process (GP) is determined by its mean function $m(x)$ and covariance function $\kappa(x,x')$. The GP can be written as $f\sim GP(m,\kappa)$. We say that a function $f$ follows a Gaussian distribution with mean function $m(x)$ and covariance function $\kappa(x,x')$. In GPR, we use the mean of the GP to approximate the unknown function. We assume that the input-output pairs are generated by a GP, which includes some parameters in the mean and covariance functions. Then, using the training data, we determine the mean and covariance parameters of the GP, obtaining the posterior GP. The posterior GP can be seen as predictors that are used to predict new inputs. One advantage of using GPR is that we can quantify the approximation function based on the covariance, thereby obtaining automatic error bounds.

We adopt the same strategy as in the previous section, except that the predictor in this section is  constructed using GPR. We denote the GPR functions by $\tilde{P}^{\omega_i}_{k,N} : P \rightarrow \tilde{\xi}^{\omega_i,\text{on}}_{k,N}$ to approximate the mapping $P^{\omega_i}_{k,N}$. After we evaluate all GPR functions $\tilde{P}^{\omega_i}_{k,N}$ for a new parameter $\mu_*$, we multiply them by the transform matrix $V$ and obtain the high fidelity eigenvectors $\tilde{\xi}^{\omega_i,\text{on}}_{k}$. 

Following the same approach, we construct the corresponding eigenfunctions $\{\tilde{\phi}_k^{\omega_i, \text{on}}: 1\leq i \leq N, 1 \leq k \leq l_i\}$ and then construct the online space $V_{\text{on}} = \text{span}\{\psi^{\omega_i,\text{on}}_k:1\leq i \leq N, 1 \leq k \leq l_i\}$, where $\psi^{\omega_i,\text{on}}_k = \chi_i\tilde{\phi}^{\omega_i,\text{on}}_k$. Finally, we obtain the coarse solution $u_H$ by solving (\ref{mumatrixpod}).

\section{The convergence of  gPC} 
\label{appendix B}
The generalized (Wiener-Askey) polynomial chaos (gPC) extends the concept of homogeneous chaos introduced by Wiener in 1938 \cite{ref25}, which originally utilized Hermite polynomials with Gaussian random variables. For any probability distribution on  $\mathbb{R}$ with finite moments of all orders, a sequence of orthonormal polynomials exists. This work focuses on continuous distributions, which are widely applicable and simplify technical complexities. We analyze chaos expansions relative to a countable set of basic random variables
 $\{\mu^m\}_n\in\mathbb{N}_0$, which need not be identically distributed but satisfy the following conditions: 
\begin{remark}\label{remark:1}
Each $\mu^m$ possesses finite moments of all orders, i.e., $E(|\mu^m|^k)<\infty$ for all $k,m\in\mathbb{N}$.
\end{remark}
\begin{remark}\label{remark:2}
The distribution functions $F_{\mu^m}(x)=P(\mu^m\le x)$  are continuous.
\end{remark}
We consider expansions in a finite set of independent random variables $\mu = (\mu^1,\mu^2,\ldots \\, \mu^M)$. 
Based on the Cameron–Martin theorem \cite{ref27} and its generalization \cite{ref28}, the gPC expansions series converge to any $L^2$ functional in the $L^2$ sense. For any random variable $\varsigma \in L^2(\Omega,\sigma(\mu),P)$ measurable with respect to $\mu$, the Doob-Dynkin lemma \cite{ref26} ensures the existence of a measurable function $f: \mathbb{R}^{\text{M}} \rightarrow \mathbb{R}$ such that $\varsigma = f(\mu)$. 
The distribution of the random M-variable $\mu$ defines a tensor measure on $\mathbb{R}^{\text{M}}$, leading to 
a probability space $(\mathbb{R},\mathcal{B}(\mathbb{R}^{\text{M}}),F_{\mu^1}(dx_1)\times F_{\mu^2}(dx_1)\ldots\times F_{\mu^M}(dx_M))$, where $\mathcal{B}(\mathbb{R}^{\text{M}})$ denotes the Borel $\sigma$-algebra on $\mathbb{R}^{\text{M}}$. This measure generates a sequence of orthonormal polynomials $\{\eta_j^{(m)}\}_{j\in\mathbb{N}_0},m=1,\ldots,M$  for each $\mu^m$, and the multivariate polynomials given by
\begin{equation}\label{eq:A.2}
      \eta_{\bm\alpha(\mu)} = \prod\limits_{m=1}^M\eta_{\alpha_m}^{(m)}(\mu^m),\ \ \ \ \ \bm\alpha =(\alpha_1,\alpha_2,\ldots,\alpha_M)\in\mathbb{N}_0^{\text{M}}. 
\end{equation}

The completeness of these polynomial systems, established by Riesz's theorem\cite{ref33}, ensures their density in $L^2$ space, provided the moment problem is uniquely solvable for each $\mu^m$. Specifically, the tensor product of orthonormal polynomial bases forms an orthonormal basis \cite{ref29} in the tensor product Hilbert space. Consequently, 
we have the following theorem:
\begin{theorem}
    Let $\mu = (\mu^1,\mu^2,\ldots,\mu^M)$ be a vector of $M \in \mathbb{N}^{\text{M}}$ independent random variables which satisfies remark A.1 and $\{\eta^{(m)}_j\}_{j\in\mathbb{N}_0},m=1,\ldots,M$, the associated orthonormal polynomial sequences. Then the orthonormal system of random variables 
    \begin{displaymath}
    \eta_{\bm\alpha(\mu)} = \prod\limits_{m=1}^M\eta_{\alpha_m}^{(m)}(\mu^m),\ \ \ \ \ \bm\alpha =(\alpha_1,\alpha_2,\ldots,\alpha_M)\in\mathbb{N}_0^{\text{M}}, 
    \end{displaymath}
    is an orthonormal basis of the space $L^2(\Omega,\sigma(\mu),P)$ if and only if the moment problem is uniquely solvable for each random variable $\mu^M$. In this case any random variable $\varsigma$ can be expanded in an abstract Fourier series of multivariate orthonormal polynomials in the basic random variables, the gPC expansion
    \begin{displaymath}
        \varsigma = \sum\limits_{\bm\alpha\in\mathbb{N}_0^{\text{M}}}a_{\bm\alpha}\eta_{\bm\alpha}{(\mu)}\ \ \ \\ \text{with coefficients} \ \ a_{\bm\alpha}= E(\varsigma\eta_{\bm\alpha}{(\mu)})
    \end{displaymath}
\end{theorem}

In summary, gPC provides a robust framework for approximating random variables using orthogonal polynomials, applicable to a broad range of probability distributions and stochastic processes. The key lies in the unique solvability of the moment problem and the completeness of the polynomial systems in their respective $L^2$ space.

\section*{Acknowledgments}
W.T. Leung is partially supported by the Hong Kong RGC Early Career Scheme 21307223. Q. Li is partially supported by National Key R \& D Program of China (No.2021YFA1001300),
National Natural Science Foundation of China (No. 12271150, 12471405).


%
 \section*{Conflict of interest}
The authors declare that they have no conflict of interest.

\bibliographystyle{unsrt}      
\bibliography{references}   

%
%

\end{document}